\documentclass[12pt,twoside]{amsart}
\usepackage{amssymb}

\nonstopmode 

\numberwithin{equation}{section}
\font\tengothic=eufm10 scaled\magstep 1
\font\sevengothic=eufm7 scaled\magstep 1
\newfam\gothicfam
      \textfont\gothicfam=\tengothic
      \scriptfont\gothicfam=\sevengothic

\newtheorem{theorem}{Theorem}[section]
\newtheorem{lemma}[theorem]{Lemma}
\newtheorem{proposition}[theorem]{Proposition}
\newtheorem{corollary}[theorem]{Corollary}
\newtheorem{conjecture}[theorem]{Conjecture}

\theoremstyle{definition}
\newtheorem{definition}[theorem]{Definition} % \theoremstyle{remark}
\newtheorem{remark}[theorem]{Remark}
\newtheorem{example}[theorem]{Example}

\newtheorem{notation}[theorem]{Notation}

\newcommand{\codim}{\operatorname{codim}}

\newcommand{\reg}{\operatorname{reg}}

\newcommand{\depth}{\operatorname{depth}}

\newcommand{\proj}[1]
{ \mathchoice
           { {\mathbb P}^{#1} }
           { {\mathbb P}^{#1} }
           { {\mathbb P}^{#1} }
           { {\mathbb P}^{#1} }
         }

\newcommand{\tor}{\operatorname{tor}}
\newcommand{\Tor}{\operatorname{Tor}}
\newcommand{\init}{\operatorname{in}}

\newcommand{\cO}{{\mathcal O}}

\newcommand {\RR}{\mathbb{R}}
\newcommand {\ZZ}{\mathbb{Z}}

\newcommand {\PP}{\mathbb{P}}

\begin{document}
\title[On the minimal free resolution of $n+1$ generic forms]{On the minimal
free resolution of $n+1$ generic forms}

\author[J.\ Migliore, R.M.\ Mir\'o-Roig]{J.\ Migliore$^*$, R.M.\
Mir\'o-Roig$^{**}$}
%\author[]{}
\address{Department of Mathematics,
        University of Notre Dame,
        Notre Dame, IN 46556,
        USA}
\email{Juan.C.Migliore.1@nd.edu}
%\author[]{U.\ Nagel}
\address{Facultat de Matem\`atiques,
Departament d'Algebra i Geometria,
Gran Via de les Corts Catalanes 585,
08007 Barcelona, SPAIN
}
\email{miro@cerber.mat.ub.es}

\date{September 20, 2001}
\thanks{$^*$ Partially supported by the University of Barcelona. \\
$^{**}$ Partially supported by DGICYT PB97-0893.}


%%%%%%%%%%%%%%%%%%%%%%%%%%%%%%%%

\begin{abstract}
Let $R = k[x_1,\dots,x_n]$ and let $I$ be the ideal of $n+1$ generically
chosen forms of degrees $d_1 \leq \dots \leq d_{n+1}$.  We give 
 the precise graded Betti numbers of $R/I$ in the following cases:
\begin{itemize}
\item $n=3$.

\item $n=4$ and $\sum_{i=1}^5 d_i$ is even.

\item $n=4$, $\sum_{i=1}^{5} d_i$ is odd and $d_2 + d_3 + d_4 < d_1 + d_5 +
4$.

\item $n$ is even and all generators have the same degree, $a$, which is even.

\item $(\sum_{i=1}^{n+1} d_i) -n$ is even and $d_2 + \dots + d_n < d_1 + 
d_{n+1} + n$.

\item $(\sum_{i=1}^{n+1} d_i) - n$ is odd, $n \geq 6 $ is even, $d_2 + \dots+d_n < d_1 +
d_{n+1} + n$ and $d_1 + \dots + d_n - d_{n+1} - n \gg 0$.
\end{itemize}
We give very good bounds on the graded Betti numbers in many other cases.  We
also extend a result of Boij by giving the graded Betti numbers for a generic
Gorenstein algebra when $n$ is even and the socle degree is large.  A recurring
theme is to examine when and why the minimal free resolution may be forced to
have redundant summands.  We conjecture that if the forms all have the same
degree then there are {\em no} redundant summands, and we present some evidence
for this conjecture.

\end{abstract}

%%%%%%%%%%%%%%%%%%%%%%%%%%%%%%%

\maketitle

\tableofcontents

%%%%%%%%%%%%%%%%%%%%%%%%%%%%%%%%%%%%%%%%%%%%%%%%%

 \section{Introduction} \label{intro}

Let $R = k[x_1,\dots,x_n]$ be a homogeneous polynomial ring over some field
$k$, and let $I = (G_1,\dots,G_d)$ be an ideal of generically chosen forms of
fixed degrees (not necessarily equal).  
A very long-standing problem in Commutative Algebra is to determine  the
Hilbert function of $R/I$.  Then a much more subtle question is to understand
all of the syzygies, i.e.\ to find the minimal free resolution of
$R/I$.  If $d \leq n$ then $I$ is a complete intersection, and its minimal free
resolution is given by the Koszul resolution.  So we assume $d > n$, which in
particular means that $R/I$ is Artinian.  A.\ Iarrobino and R.\ Fr\"oberg have
made conjectures about the Hilbert function, and A.\ Iarrobino has made a 
conjecture for the minimal free resolution in this case.  One of the
consequences of our work is to give a counterexample to the latter conjecture.

Several contributions to this very difficult problem have been made.  We
first discuss the Hilbert function.  If $d = n+1$ then the Hilbert function
is well known, coming from a result of R.\ Stanley \cite{stanley} and of J.\
Watanabe \cite{watanabe} which implies that a general Artinian complete
intersection has the Strong Lefschetz property (cf.\ Definition \ref{def of
wlp and slp}).  We will discuss this  shortly, but for now we
do not yet assume $d = n+1$.  The case
$n=2$ was solved by Fr\"oberg \cite{froberg}.  The case $n=3$ was solved by
Anick \cite{anick}.  M.\ Hochster and D.\ Laksov \cite{hochster-laksov}
showed that a generically chosen set of forms of the same degree span as
much as possible in the next degree.  (Note that this gives the value of the
Hilbert function in the next degree, and it also gives the number of linear
syzygies of the forms.)  This was extended by Aubry \cite{aubry}.  Also,
Fr\"oberg and Hollman \cite{froberg-hollman} solved it for forms of degree 2
if $n \leq 11$, and for forms of degree 3 if $n \leq 8$.

Apart from the above work, nothing seems to be published about the problem
of finding the minimal free resolution for generic forms.  This is the
central problem which we address in this paper.  We remark that very
different approaches to this subject are being carried out by Ben Richert
and Keith Pardue  \cite{PR} and by Karen Chandler \cite{C}.

A related problem is the Minimal Resolution Conjecture \cite{lorenzini}.  A
generic set of points in projective space has so-called {\em generic Hilbert
function}, which depends only on the number of points.  The Minimal
Resolution Conjecture asked whether the entire resolution similarly is the
``expected'' one, in the sense that the graded Betti numbers depend only on
the Hilbert function.  In particular, it requires that there not be any
``ghost'' terms in the resolution, i.e.\ that consecutive terms in the
resolution never have a summand in common.  (Ghost terms cannot be detected
from the Hilbert function alone.)  Unfortunately, it was shown by D.\
Eisenbud and S.\ Popescu \cite{EP} that this is not true.  The first
counterexample is the case of 11 points in $\proj{6}$, discovered in
computational experiments by F.\ Schreyer, where there is a summand $R(-5)$
in both the third and fourth syzygies which does not split off.  

Eleven general points in $\proj{6}$ lie on 17 independent quadrics.  A
natural question is then whether 17 general quadrics in $k[x_1,\dots,x_6]$
also pick up a ghost term, or whether in fact it has the expected
resolution.  One can check that in fact it does have the expected
resolution (cf.\ \cite{IK} page 197).  So the Artinian reduction of 11
general points in $\proj{6}$ is not ``general'' enough as an Artinian
algebra.  

This question of whether ghost terms exist in the minimal free
resolution was of central interest to us in writing this paper.  It is clear
that they cannot be entirely avoided.  For instance, if our chosen degrees
include two forms of degree 4 and one of degree 8, then we naturally expect
a Koszul syzygy of degree 8, so there is a summand $R(-8)$ in the first
syzygy module which does not split off with the summand $R(-8)$
corresponding to the generator.  

A natural conjecture, due to Iarrobino \cite{iarrobino} , is that the ghost
terms arising as a result of Koszul syzygies should be the only kind of
exception.  Called the Thin Resolution Conjecture, it says that  ``the minimal
free resolution
\dots is the minimum one that is consistent with their (expected) Hilbert
function; that is, the Koszul resolution up to the smallest degree where''
$(R/I)_i = 0$ (cf.\ \cite{IK} page 197).  

One result of our work is a clearer understanding of the fact that other
ghost terms do in fact arise!  For instance, we show in Example \ref{codim 3
ghosts} that when $n=3$ and generic forms of degree 4, 4, 4 and 8 respectively,
the minimal free resolution is
\[
0 \rightarrow
\left (
\begin{array}{c}
R(-10) \\
\oplus \\
R(-11)^2
\end{array}
\right )
\rightarrow
\left (
\begin{array}{c}
R(-8)^3 \\
\oplus \\
R(-9)^2 \\
\oplus \\
R(-10)
\end{array}
\right )
\rightarrow
\left (
\begin{array}{c}
R(-4)^3 \\
\oplus \\
R(-8)
\end{array}
\right )
\rightarrow R \rightarrow R/I \rightarrow 0
\]
We see that the $R(-8)$ does not split, as predicted above, but that
furthermore there is a summand $R(-10)$ shared by the second and third
modules which also does not split.  Notice that the Hilbert function of $R/I$ is
\[
1  \ \ 3 \ \ 6 \ \ 10 \ \ 12 \ \ 12 \ \ 10 \ \ 6 \ \ 2,
\]
and that the summand $R(-10)$ does not correspond to a Koszul syzygy.

Iarrobino informs us that the above example is a counterexample to his Thin
Resolution Conjecture, and furthermore that it is a counterexample to his
published statement \cite{iarrobino} that the Thin Resolution Conjecture had
been shown to be equivalent to Fr\"oberg's Conjecture on the Hilbert function.
Other examples of ghost terms that arise can be found in Example \ref{no
splitting}, Example \ref{codim 3 ghosts} and Example \ref{3 3 4 6 6}, but
one can produce more from the theorems.

This paper concerns solely the case of $n+1$ generic forms in
$k[x_1,\dots,x_n]$, i.e.\ an Artinian almost complete intersection.  We let
$I = (G_1,\dots, G_{n+1})$ where $\deg G_i = d_i$ and $d_1 \leq \dots \leq
d_{n+1}$.  We always assume that $d_{n+1} \leq (\sum_{i=1}^n d_i) -n$ because
otherwise $G_{n+1}$ is in the ideal generated by the first $n$ generators,
and so $I$ is a complete intersection.

Our
first observation (which is not new) is that such an ideal $I$ can be linked
to a Gorenstein ideal $G$ via the complete intersection $J$ defined by the
first $n$ generators of $I$.  In Lemma \ref{G hilb function facts} we give
some facts about the Hilbert function of $R/G$.  One is to note that the
Hilbert function of $R/G$ has either one ``peak'' in the middle, or two. The
technical condition for one peak is that $(\sum_{i=1}^{n+1} d_i ) -n$ must be
even.  Furthermore, we describe exactly when $R/G$ agrees with $R$ all the
way up to this peak.  For short, we will say that ``the growth of
$R/G$ is maximal'' in this case.  The technical condition for the growth to
be maximal is $d_2+\dots+d_n < d_1 + d_{n+1} + n$.

As a consequence of Lemma \ref{G hilb function facts}, we show
in Corollary \ref{wlp} that $A := R/G$ has the Strong Lefschetz
Property.  This is central
especially for Section 5, where we have our strongest results, because it
allows us to compute the Hilbert function of $A/LA$ for a general linear
form $L$.

We observe in Section 3 that a free resolution for $R/I$ can be given in
terms of one for $R/G$ (again this is not new), and that we can control to a
large degree the possible splitting.  So the problem is reduced to finding a
minimal free resolution for $R/G$.  Sections 3, 4 and 5 give different
approaches to this, for different situations.

In Section 3 we first use a result of the first author and U.\ Nagel
\cite{MN3}, which gives the precise minimal free resolution for $R/G$ when
its Hilbert function has only one peak and the growth is maximal.  We then
determine exactly what splitting can occur for the linked ideal, giving the
minimal free resolution for $R/I$ (Corollary \ref{one peak resol}). 

The more difficult situation (still assuming that the growth is maximal) is
when the Hilbert function of $R/G$ has two peaks.  Here the results of
\cite{MN3} do not give sharp bounds on the graded Betti numbers for $R/G$. 
However, a result of Boij \cite{boij} on generic Gorenstein ideals is helpful
here when $n=4$.  We generalize Boij's result, giving the minimal free
resolution of a generic Gorenstein algebra when $n$ is even and the socle
degree is large (Proposition \ref{gor two peaks even}).  (We give bounds
when $n$ is odd, in Remark \ref{gor two peaks odd}.)  As a consequence we
give the precise minimal free resolution for $R/I$ (Proposition \ref{aci two
peaks}) when the Hilbert function of $R/G$ has two peaks, the growth is
maximal, $n$ is even and the socle degree is large.

In section 4 we give a complete answer to the resolution problem for $n=3$. 
Our method is to apply the work of Diesel \cite{diesel} to find the minimal
free resolution for the generic Gorenstein algebra with the known Hilbert
function (coming only from the choice of $d_1,\dots,d_4$), and then apply
our methods to determine all the splitting that occurs.  The main result here
is  Theorem \ref{resol for codim 3}.

In Section 5 we use a different result of the first author and Nagel
\cite{MN3} to make a more subtle study of the minimal free resolution of
$A = R/G$.  The procedure is the following.  First determine the Hilbert
function of $A/LA$ for a general linear form $L$, which is known from the
Weak Lefschetz property.  Then determine the graded Betti numbers of $A/LA$
over $R/(L)$.  This information, together with the result from \cite{MN3},
allows us to make very good bounds for the graded Betti numbers of $A$.  A
careful analysis then shows that these bounds are actually sharp!!  Finally,
the link to $I$ is studied, and it is determined exactly what splitting
occurs, resulting in the minimal free resolution of $R/I$.  This program
gives the following (Theorems \ref{main result of section 5} and
\ref{general ans when n=4}):

\begin{itemize}
\item Assume all $n+1$ generators have the same degree, $a$.  Let $s(n,a)
= (n-1)a -n$.  Then
\begin{itemize}
\item If $n$ is odd, we give a resolution for $R/I$ that is not quite
minimal.
\item If $n$ is even and $s(n,a)$ is odd, we give a resolution for $R/I$
that is not quite minimal.
\item If $n$ is even and $s(n,a)$ is even then we give the precise minimal
free resolution for $R/I$.
\end{itemize}

\item Assume that $n=4$ and that $\sum_{i=1}^5 d_i$ is even.  Then we give
the precise minimal free resolution for $R/I$.
\end{itemize}

We note that the approach of this section can be applied in other
situations, but that the notation quickly becomes overwhelming.

We would also like to remark that all of our results also hold in another
context.  Instead of general forms of degree $d_i$, fix generally chosen
linear forms $L_1,\dots, L_{n+1}$ and consider $G_i = L_i^{d_i}$.  Then it
is still true that the ideal $J = (G_1,\dots, G_n)$ has the Strong Lefschetz
property (this is the original result of Stanley and of Watanabe).  Then all
of the machinery of this paper carries through to this setting.  This gives
a connection to the study of fat points, for which we refer to work of
Chandler and of Iarrobino.

Most of this work was done while the first author was a guest of the
University of Barcelona, and he would like to thank the professors and
students of the Departament d'Algebra i Geometria for their warm
hospitality.  The authors also thank Karen Chandler, Tony Geramita and Tony
Iarrobino for helpful references and comments.  

\section{Hilbert function calculations}

Let $R = k[x_1,\dots,x_n]$ where $k$ is an algebraically closed field
(although we remark that this hypothesis is needed only for \S 3 beginning
with Proposition \ref{gor two peaks even}).   For any homogeneous ideal
$I \subset R$ we denote the Hilbert function of $R/I$ by
$h_{R/I}(t)$.  If $R/I$ is Gorenstein, we sometimes refer to $I$ itself as
being Gorenstein.  In this paper, for a numerical function $f$ we denote by
$\Delta f$ the first difference function $\Delta f(t) = f(t) - f(t-1)$ for
all $t \in {\mathbb Z}$.  

\begin{definition}\label{basic-Gor-defs}
 Let $\underline{h} = (h_0,\dots,h_s)$ be a sequence of positive integers.
$\underline{h}$ is called a {\em Gorenstein sequence} if it is the Hilbert
function of some Gorenstein Artinian $k$-algebra.  $\underline{h}$ is {\em
unimodal} if $h_0 \leq h_1 \leq \dots \leq h_j \geq h_{j+1} \geq \dots \geq
h_s$ for some $j$.  $\underline{h}$ is called an {\em SI-sequence} (for
Stanley-Iarrobino) if it satisfies the following two conditions:

\begin{itemize}
\item[(i)] $\underline{h}$ is symmetric, i.e.\ $h_{s-i} = h_i$ for all
$i=0,\dots,\lfloor \frac{s}{2} \rfloor$.
\item[(ii)] $(h_0, h_1-h_0,h_2-h_1,\dots,h_j-h_{j-1})$ is an O-sequence,
where $j = \lfloor \frac{s}{2} \rfloor$; i.e.\ the ``first half'' of
$\underline{h}$ is a {\em differentiable} O-sequence.
\end{itemize}
\end{definition}

\begin{definition} \label{def of wlp and slp}
A standard graded Artinian $k$-algebra $A = \bigoplus_{i \geq 0} A_i$ has the
{\em Weak Lefschetz property} (sometimes called the {\em Weak Stanley
property}) if for each $i$, the multiplication $A_i \rightarrow A_{i+1}$
induced by a general linear form $L$ has maximal rank.  $A$ has the {\em
Strong Lefschetz property} if for each $i$ and each $d \geq 1$, the
multiplication $A_i \rightarrow A_{i+d}$ induced by a $L^d$, for a general
linear form $L$, has maximal rank.
\end{definition}

It was shown by Harima \cite{harima} that the SI-sequences characterize the
possible Hilbert functions of graded Artinian Gorenstein $k$-algebras with
the Weak Lefschetz property.  The set of all possible Hilbert functions for
Artinian $k$-algebras with the Weak or Strong Lefschetz properties was
described in \cite{HMNW}.

\begin{notation}
Let $I = (G_1,\dots,G_n, G_{n+1}) \subset R$ be an ideal where
$\deg G_i = d_i$ for $1 \leq i \leq n+1$ and for each $i$, $G_i$ is chosen
generically in $R_{d_i}$.  We will call $I$ a {\em general Artinian almost
complete intersection of type $(d_1,\dots,d_n,d_{n+1})$}.
\end{notation}

\begin{remark} \label{hyp}
Let $I$ be a general Artinian almost complete intersection of type \linebreak
$(d_1,\dots,d_n,d_{n+1})$.  Without loss of generality assume that $d_1 \leq
\dots \leq d_{n+1}$. Because the forms are chosen generically, we may assume
that $J := (G_1,\dots,G_n)$ forms a regular sequence.  If $d_{n+1} >
(\sum_{i=1}^n d_i) - n$ then $G_{n+1} \in J$ (since it is Artinian) and hence
$I = J$, and the Hilbert function and minimal free resolution of $I$ are hence
known (from the Koszul resolution).  So without loss of generality, {\em from
now on we assume that $d_{n+1} \leq (\sum_{i=1}^n d_i) - n$.  }
\end{remark}

The Hilbert function of $R/I$ is well known, and we remind the reader of the
main idea.  Stanley \cite{stanley} and Watanabe \cite{watanabe} independently
showed that a general Artinian complete intersection, $J$, has the Strong
Lefschetz property.  Hence we have
$h_{R/I}(t) = \max \{ h_{R/J}(t) - h_{R/J}(t-d_{n+1}),0\}$, and
the values of $h_{R/J}(t)$ are known thanks to the Koszul resolution. In
particular, we have:

\begin{lemma} \label{end of R/I}
The socle degree of $R/I$ (i.e.\ the degree of the last non-zero component of
$R/I$) is
\[
\left \lfloor \frac{1}{2} \left ( \left ( \sum_{i=1}^{n+1} d_i  \right ) -n
-1 \right )
\right
\rfloor
\]
where $\lfloor t \rfloor$ denotes the greatest integer less than or equal to
$t$.
\end{lemma}

\begin{proof}
We have assumed that $d_1 \leq \dots \leq d_n \leq d_{n+1}$.  Furthermore, we
know that the Hilbert function of $R/J$ is symmetric and ends in degree
$\sum_{i=1}^n d_i - n$.  Hence the ``flat part'' of the Hilbert function has
length less than $d_{n+1}$.  Because of the Strong Lefschetz property, the
degree we are looking for is the greatest $j$ for which $\dim (R/J)_j > \dim
(R/J)_{j-d_{n+1}}$.  This is a slightly tedious but easy computation, checking
different parity cases.  It requires only the facts mentioned and not the
precise values of these dimensions.
\end{proof}

The ideal $J$ links $I$ to an ideal $G := [J : I]$ which is easily seen to be
arithmetically Gorenstein, using the standard mapping cone construction (cf.\
\cite{PS}, \cite{migliore}).  Since the Hilbert function of $R/J$ and of
$R/I$ are known, we can compute the Hilbert function of $R/G$ (cf.\
\cite{DGO}, \cite{migliore}).  More precisely, we have

\begin{equation}  \label{hilbftnform}
\begin{array}{rcl}
h_{R/G}(t) & = &  h_{R/J}\left ( (\sum_{i=1}^n d_i) - n - t \right
) - h_{R/I} \left ( (\sum_{i=1}^n d_i) - n - t \right ) \\ \\
& = & h_{R/J}\left ( (\sum_{i=1}^n d_i) - n - t \right
)  \\ \\
&& \hbox{\hskip .5cm} - \left [ h_{R/J} \left ( (\sum_{i=1}^n d_i) - n - t
\right ) - h_{R/J}
\left ( (\sum_{i=1}^n d_i) - n - t - d_{n+1} \right ) \right ]_+
\end{array}
\end{equation}
where we denote by $[x]_+$ the maximum of $x$ and zero.
Since $R/G$ is Gorenstein, this Hilbert
function is symmetric, so we only have to compute half of it.

For our purposes below we need to know precisely for how many degrees $t$
does the Hilbert function $h_{R/G}(t)$ achieve its maximum value.  We
will say that the Hilbert function {\em has $r$ peaks} if there are $r$ such
degrees, and we will note that these must occur consecutively (i.e. the
Hilbert function is unimodal).  We also will need to know under what
conditions the growth of the Hilbert function is the maximum possible (i.e.\
coincides with the polynomial ring) for the entire first half of the Hilbert
function.

\begin{lemma} \label{G hilb function facts}
Let $G$ be the arithmetically Gorenstein ideal linked to $I$ by the complete
intersection $J$ as above.
\begin{itemize}
\item[(a)] The socle degree of $R/G$ is
$(\sum_{i=1}^n d_i) - d_{n+1} - n.$

\item[(b)] The Hilbert function $h_{R/G}$ is unimodal, with one peak if
$(\sum_{i=1}^{n+1} d_i ) -n$ is even, and two peaks if it is odd.

\item[(c)] For all integers $t \leq \frac{(\sum_{i=1}^n d_i) - d_{n+1} -
n}{2}$ we have $h_{R/G}(t) = h_{R/J}(t)$.  By symmetry of $h_{R/G}$, this
completely determines $h_{R/G}$.

\item[(d)] We have
\[
h_{R/G}(t) = \binom{t+n-1}{n-1} \hbox{ for all integers } 0 \leq t \leq
\frac{(\sum_{i=1}^n d_i) - d_{n+1} - n}{2}
\]
if and only if $d_2+\dots+d_n < d_1 + d_{n+1} + n$.
\end{itemize}
\end{lemma}

\begin{proof}
Clearly $R/I$ and $R/J$ first differ in degree $d_{n+1}$.  Since $R/J$ ends
in degree \linebreak $(\sum_{i=1}^n d_i) - n$, (a) follows immediately from a
Hilbert function calculation as indicated in (\ref{hilbftnform}).

For (b), note first that the parity does not change if we replace
$\sum_{i=1}^{n+1} d_i$ by \linebreak $\sum_{i=1}^n d_i - d_{n+1}$.  Because
of (a) and the symmetry of the Hilbert function of $R/G$, we know that the
number of peaks of $R/G$ will be odd if $(\sum_{i=1}^{n+1} d_i ) -n$ is even,
and even if $(\sum_{i=1}^{n+1} d_i ) -n$ is odd (once we have shown
unimodality).

First assume that $(\sum_{i=1}^{n+1} d_i) -n$ is even.  We want to show that
$h_{R/G}$ is unimodal with one peak, which by (a) and symmetry would have
to occur in degree $\frac{(\sum_{i=1}^n d_i) - d_{n+1} - n}{2}$.  Using Lemma
\ref{end of R/I} and the formula (\ref{hilbftnform}), one quickly can check
that for any $t \geq 0$ we have
\[
h_{R/G}\left ( \frac{(\sum_{i=1}^n d_i) - d_{n+1} - n}{2} -t \right ) =
h_{R/J}\left ( \frac{(\sum_{i=1}^{n+1} d_i) -n}{2} +t \right ).
\]
Because of our hypothesis that $d_{n+1} \leq (\sum_{i=1}^n d_i) -n$, one can
check that $d_n \leq d_{n+1} \leq \frac{(\sum_{i=1}^{n+1} d_i) -n}{2}$.  Hence
for $t \geq 0$ the right-hand side of the above equation is strictly
decreasing.  This proves (b) for $(\sum_{i=1}^n d_i) - d_{n+1} - n$ even.
Furthermore, (c) is easy to check using the symmetry of $h_{R/J}$.

If $(\sum_{i=1}^{n+1} d_i) -n$ is odd, note that the fraction in the
statement of (c) is not an integer.  The proof of (b) and (c) is identical to
that of the previous case, simply replacing $\frac{(\sum_{i=1}^n d_i) -
d_{n+1} - n}{2}$ by $\frac{(\sum_{i=1}^n d_i) - d_{n+1} - n-1}{2}$.

The proof of (d) follows immediately from (c) and a calculation, by setting
$d_1 > \frac{\sum_{i=1}^n d_i - d_{n+1} - n}{2}$ and simplifying.
\end{proof}

\begin{corollary} \label{wlp}
$R/G$ has the Strong Lefschetz property. 
\end{corollary}

\begin{proof} 
We have $J = (G_1,\dots,G_n)$, $I = J+(G_{n+1})$ and $G = [J:I]$ in the
polynomial ring $R = K[x_1,\dots,x_n]$.  Thanks to the result of Stanley
\cite{stanley} and of Watanabe
\cite{watanabe},
$R/J$ has the Strong Lefschetz property.  Now consider the exact sequence
(not graded)
\[
\begin{array}{rcl}
0 \rightarrow [0:_{R/J} G_{n+1}] \rightarrow 
& R/J  \  \stackrel{G_{n+1}}{\longrightarrow} \  R/J &
\rightarrow R/I \rightarrow 0 \\
& \hbox{\hskip .5cm} \searrow \hfill \nearrow \hbox{\hskip .5cm}  \\
& I/J  \\ 
& \hbox{\hskip .5cm}  \nearrow \hfill \searrow  \hbox{\hskip .5cm}  \\
& 0 \hbox{\hskip 2.3cm} 0
\end{array}
\]
By Theorem 4.14 of \cite{IK}, if $A$ is a graded algebra with the Strong
Lefschetz property then so is $A/[0:F]$, where $F$ is a general element of
$A$.  Therefore, the cokernel $I/J$ has the Strong Lefschetz property.  But
$J$ is a complete intersection linking $I$ to $G$, so we have an isomorphism
\[
I/J \cong K_G (n-d)
\]
where $K_G$ is the canonical module of $R/G$ and $d = \sum_{i=1}^n d_i$. 
Since $R/G$ is Gorenstein, $K_G$ is isomorphic to $R/G$ up to twist.  It
follows that $R/G$ has the Strong Lefschetz property as claimed.
\end{proof}

%%%%%%%%%%%%%%%%%%%%%%%%%%%%%%%%%%%%%%%%%%%%%%%%%%%%%%%%%%%%%%%%%%%%%%%%%%%%%

\section{First Bounds for the graded Betti numbers of an almost complete
intersection}

In this section we are interested in describing the minimal free resolution
for a general almost complete intersection in $R = k[x_1,\dots,x_n]$.

First we consider the case where one of the generators has degree 1.  Our
approach is the same as that of Lemma 2.6 and
Corollary 2.7 of \cite{HMNW}, although the result we obtain in (b) is not
explicit, as it is in \cite{HMNW} for the case of three variables.  We give
the precise graded Betti numbers if
$d_{n+1}$ is large enough, and we show how to reduce it to the same problem in
a smaller polynomial ring otherwise.

\begin{proposition} \label{one gen of deg 1}
Let $I = (L,G_2,\dots,G_{n+1})$ be a general almost complete intersection in
$R = k[x_1,\dots,x_n]$, with $d_1 = \deg L = 1$, $d_i = \deg G_i$ for $i \geq
2$, and $d_2 \leq \dots \leq d_n \leq d_{n+1}$.  Let $J =
(G_2,\dots,G_{n+1})$,  $\bar R = R/(L)$ and let $\bar J \subset \bar R$ be the
reduction of $J$ modulo $L$.  Note that $\bar J$ is a general almost complete
intersection in $\bar R$ of type $(d_2,\dots,d_{n+1})$.

\begin{itemize}
\item[(a)] If  $d_{n+1} > (\sum_{i=1}^n d_i) -n$ then  $G_{n+1} \in
(L,G_2,\dots,G_n)$, so $I = (L,G_2,\dots,G_n) \subset R$ is an Artinian
complete intersection, and its minimal free resolution over $R$ is given by
the Koszul resolution.

\item[(b)] If $d_{n+1} \leq (\sum_{i=1}^n d_i) -n$, let $J' = R \bar J
\subset R$.  Then $J'$ is the saturated ideal of an almost complete
intersection ideal in $R$ with $\depth R/J' = 1$, and the minimal free
resolution of $R/I$ is given by the tensor product of $R/J'$ and $R/(L)$.  In
particular,
\[
[\tor_i^R (R/I,k)]_j = [\tor_i^{\bar R} (\bar R/\bar J, k)]_j +
[\tor_{i-1}^{\bar R} (\bar R/ \bar J,k)]_{j+1}
\]
for $1 \leq i \leq n$.
\end{itemize}
\end{proposition}

\begin{proof}
Part (a) is immediate from Remark \ref{hyp}.  For part (b), note that $I =
J+ (L) = J' + (L)$.  The fact that $J'$ is Cohen-Macaulay of depth one follows
since $\bar J$ is Artinian in $\bar R$, so $L$ is a non-zero divisor for
$R/J'$.  So the graded Betti numbers of $J'$ over $R$ are the same as those of
$\bar J$ over $\bar R$.  Note that $\Tor_i^{\bar R} (\bar R/\bar J,k) =
\Tor_i^R (R/J',k) = 0$ for $i \geq n$.
\end{proof}

For the remainder of this section we will strongly use some results from
\cite{MN3}.  We first collect some notation.

\begin{notation} \label{ring notation}
Let $R = k[x_1,\dots,x_n]$ and let $I \subset R$ be a general Artinian almost
complete intersection of type $(d_1,\dots,d_n,d_{n+1})$.
  As before, let $J$ be the complete intersection given by the first $n$
generators of $I$, and let $G$ be the linked Gorenstein ideal. Let $A = R/G$.

 For a graded $R$-module $M$, set $[\tor_i^R (M,k)]_j := \dim_k
[\Tor_i^R(M,k)]_j$.

We set
\[
\begin{array}{rcl}
s  &  = & \displaystyle  \sum_{i=1}^n d_i - d_{n+1} - n,       \\  \\  \alpha
& = & \displaystyle \left \lfloor \frac{s+1}{2}  \right \rfloor \end{array}
\] Note
that $s$ is the socle degree of $A$ (Lemma \ref{G hilb function
facts}) and $\alpha $ is the degree in which the last peak (of
which there are either one or two) occurs in the Hilbert function
of $A$. Since $A$ has the Weak Lefschetz property (even the Strong
Lefschetz property -- see Corollary
\ref{wlp}) then we also have $\alpha=in[0:_A L]$ for a general linear form
$L$, where for a homogeneous ideal $I$, $in(I)$ is its initial degree. 
%
% If $A$ does not
%have the Weak Lefschetz property, then by Corollary \ref{wlp} we have
%$\alpha-1 = in[0:_A L]$.
\end{notation}

We will need the following results, which we specialize slightly for our
purposes.

\begin{notation}
Let $\underline{h} = (1,h_1,\dots,h_s)$ be the $h$-vector of an Artinian
$k$-algebra. Let $c \geq h_1$ be an integer.  Then there is a uniquely
determined lex-segment ideal $I \subset k[z_1,\dots,z_{c}] =: T$ such that
$T/I$ has
$\underline{h}$ as its Hilbert function.  We define
\[
\beta_{i,j} (\underline{h},c) := \left [ \tor_i^T (T/I,k) \right ]_j.
\]
If $c=h_1$ we simply write $\beta_{i,j}(\underline{h})$ instead of
$\beta_{i,j}(\underline{h},h_1)$.
\end{notation}

Using the main result of \cite{bigatti} and \cite{hulett}, the following was
obtained:

\begin{theorem}[\cite{MN3} Theorem 8.13] \label{MN3 main result}
Let $\underline{h} = (1,h_1,h_2,\dots,h_\ell,\dots, h_s)$ be an SI-sequence
where
$h_{\ell-1 } < h_\ell = \cdots = h_{s-\ell} > h_{s-\ell+1}$.  Put
$\underline{g} = (1,h_1-1,h_2-h_1,\dots,h_\ell-h_{\ell-1})$.  Then we have
\begin{itemize}
\item[(a)] If $A = R/G$ is a Gorenstein $k$-algebra with $c = \codim G$ and
having  $\underline{h}$ as $h$-vector and an Artinian reduction
which has the Weak Lefschetz property, then
\[
\left [ \tor_i^R (A,k) \right ]_j \leq
\left \{
\begin{array}{ll}
\beta_{i,j} (\underline{g},c-1) & \hbox{if $j \leq s-\ell+i-1$} \\
\beta_{i,j} (\underline{g},c-1) + \beta_{c-i,s+c-j} (\underline{g},c-1) &
\hbox{if $s-\ell+i \leq j \leq \ell+i$} \\
\beta_{c-i,s+c-j} (\underline{g},c-1) & \hbox{if $j \geq \ell+i+1$}
\end{array}
\right.
\]
\item[(b)] Suppose that the field $k$ has ``sufficiently many elements''
(e.g.\ infinitely many).  Then there is a reduced, non-degenerate
arithmetically Gorenstein subscheme $X \subset \proj{n} = \hbox{Proj} (R)$ of
codimension $c$ whose Artinian reduction has the Weak Lefschetz property
and $h$-vector $\underline{h}$, and with equality holding for all $i,j \in
{\mathbb Z}$ in {\rm (a)}.  This subscheme can be constructed explicitly.
\end{itemize}
\end{theorem}

\begin{corollary}[\cite{MN3} Corollary 8.14] \label{precise resol}
Let $s,\ell$ be positive integers, where either $s = 2\ell$ or $s \geq
2\ell+2$.  Define $\underline{h} = (h_0,\dots,h_s)$ by
\[
h_i = \left \{
\begin{array}{rl}
\binom{n-1+i}{n-1} & \hbox{if $0 \leq i \leq \ell$}; \\
\binom{n-1+\ell}{n-1} & \hbox{if $\ell \leq i \leq s-\ell$};\\
\binom{s-i+n-1}{n-1} & \hbox{if $s-\ell \leq i \leq s$}
\end{array}
\right.
\]
Let $G \subset R$ be an Artinian Gorenstein ideal such that $R/G$
has the Weak Lefschetz property and Hilbert function
$\underline{h}$. Then the minimal free $R$-resolution of $R/G$ has
the shape
\[
0 \rightarrow R(-s-n) \rightarrow F_{n-1} \rightarrow \cdots \rightarrow F_1
\rightarrow R \rightarrow R/G \rightarrow 0
\]
where
\[
F_i = R(-\ell-i)^{\alpha_i} \oplus R(\ell-s-i)^{\gamma_i} \ \ \hbox{ and } \ \
\alpha_i = \binom{n+\ell-1}{i+\ell} \binom{\ell-1+i}{\ell} = \gamma_{n-i}.
\]
\end{corollary}

\begin{remark} \label{observations about MN3}
\begin{itemize}
\item[a.] In the language of this paper, $\ell$ refers to the first degree
where the peak is achieved, $s$ is the last degree in which the Hilbert
function is non-zero, and $s-\ell$ is the last degree where there is a peak.
The condition
$s=2\ell$ means that there is one peak, while the case $s \geq 2\ell+2$ means
three or more peaks (which does not occur for us).  In our situation,
\[
\begin{array}{rcl}
s & = & (\sum_{i=1}^n d_i )- d_{n+1} - n \\ \\
\ell & = & \left \{
\begin{array}{ll}
\frac{(\sum_{i=1}^n d_i) - d_{n+1} - n}{2} & \hbox{ if $(\sum_{i=1}^{n+1} d_i)
 - n$ is even;}\\ \\
\frac{(\sum_{i=1}^n d_i) - d_{n+1} - n-1}{2} & \hbox{ if $(\sum_{i=1}^{n+1}
d_i)
 - n$ is odd.}
\end{array}
\right.
\end{array}
\]

\item[b.] The Hilbert function described in this corollary is the same as that
of Lemma \ref{G hilb function facts} (d).

\item[c.] The hypothesis given in Corollary \ref{precise resol} that $R/G$
have the Weak Lefschetz property is not needed in the case $s=2\ell$, since
the growth described is precisely that of $R$ up to degree $\ell$: this means
that multiplication by a general linear form is injective in this range, and
the surjectivity in the other degrees comes from the self-duality of $R/G$.
In any case, we have seen that our $G$ even has the Strong Lefschetz
property, hence in particular the Weak Lefschetz property.
\end{itemize}
\end{remark}

We now make the basic construction of this section, and then we
will draw consequences.  We have a general almost complete
intersection $I = (G_1,\dots,G_n,G_{n+1})$ with $d_1 \leq \dots
\leq d_n \leq d_{n+1} \leq (\sum_{i=1}^n d_i) - n$, the complete
intersection $J = (G_1,\dots,G_n) \subset I$ and the Gorenstein
ideal $G = [J:I]$.  Let $d = d_1 + \dots + d_n$.  Consider the
minimal free $R$-resolution of $R/J$ given by the Koszul
resolution:
\[
0 \rightarrow K_n \rightarrow K_{n-1} \rightarrow \cdots \rightarrow K_2
\rightarrow K_1 \rightarrow R \rightarrow R/J \rightarrow 0
\]
where
\[
K_i = \bigwedge^i \left ( \bigoplus_{i=1}^n R(-d_i) \right ).
\]
In particular, $K_1 \cong \bigoplus_{i=1}^n R(-d_i)$, $K_{n-1} \cong
\bigoplus_{j=1}^n R(d_j - d)$ and $K_n \cong R(-d)$.  Consider the minimal
free resolution of $R/G$:
\[
0 \rightarrow R(-e) \rightarrow F_{n-1} \rightarrow \cdots \rightarrow F_2
\rightarrow F_1 \rightarrow R \rightarrow R/G \rightarrow 0
\]
where $e = (\sum_{i=1}^n d_i) - d_{n+1}$ thanks to Lemma \ref{G hilb function
facts} (a).

Applying the mapping cone construction to the diagram
\[
\begin{array}{ccccccccccccccccccccc}
0 & \rightarrow & R(-d) & \rightarrow & K_{n-1} & \rightarrow & \cdots &
\rightarrow & K_2 & \rightarrow & K_1 & \rightarrow & R
& \rightarrow & R/J & \rightarrow & 0 \\
&& \downarrow && \downarrow &&&& \downarrow && \downarrow && \downarrow &&
\downarrow \\
0 & \rightarrow & R(-e) & \rightarrow & F_{n-1} & \rightarrow & \cdots &
\rightarrow & F_2 & \rightarrow & F_1 & \rightarrow & R
& \rightarrow & R/G & \rightarrow & 0
\end{array}
\]
gives a free $R$-resolution for $R/I$:
\begin{equation} \label{resol of R/I}
0 \rightarrow F_1^\vee (-d) \rightarrow
\begin{array}{c}
K_1^\vee (-d) \\
\oplus \\
F_2^\vee (-d)
\end{array}
\rightarrow
\begin{array}{c}
K_2^\vee (-d) \\
\oplus \\
F_3^\vee (-d)
\end{array}
\rightarrow \cdots \rightarrow
\begin{array}{c}
K_{n-1}^\vee (-d) \\
\oplus \\
R(e-d)
\end{array}
\rightarrow R \rightarrow R/I \rightarrow 0
\end{equation}
A simple calculation gives that $I$ has $n+1$ generators, in degrees
$d_1,\dots,d_n,d_{n+1}$, as expected.  The challenge is to determine how much
splitting can occur and to try to narrow down as much as possible what the
free modules $F_i$ can be.  In any case, we immediately obtain

\begin{proposition} \label{first bound}
With the above notation, we have
\[
[\tor_i^R (R/I,k)]_j \leq [\tor_{n-i}^R (R/J,k)]_{d-j} + [\tor_{n-i+1}^R
(R/G,k)]_{d-j}.
\]
\end{proposition}

The first term on the right-hand side of the inequality in Proposition
\ref{first bound} is obtained from the Koszul resolution and is completely
determined.  The second term can be bounded as follows:

\begin{corollary} \label{weak bound} With the above notation we have
\[
\tor_i^R (R/I,k)_j \leq \tor_{n-i}^R (R/J,k)_{d-j} + B_{n-i+1,d-j}
\]
where $B_{a,b}$ is the bound for $[\tor_a^R (R/G,k)]_b$ obtained in Theorem
\ref{MN3 main result}.
\end{corollary}

Corollary \ref{weak bound} is slightly too general for our purposes (usually).
It assumes that the ``first half'' of $R/G$ is as bad as possible, while in
our case it agrees with a complete intersection.
The ``obvious'' first place to look for splitting in Proposition
\ref{first bound} is to see if any of the generators $G_i$ of $J$ are minimal
generators of
$G$.  If this occurs, each such generator leads to a splitting of a rank one
free summand at the end of the resolution of $R/I$.
This need not happen, however, as illustrated by the following example.

\begin{example} \label{ex of no split}
Let $n=4$, $d_1= \dots = d_4 = 5$, $d_5 = 10$.  By Lemma \ref{G hilb function
facts} (b) this has one peak, which one computes occurs in degree  $\ell =
3$.  By Lemma \ref{G hilb function facts} (d) we have that the Hilbert
function of
$R/G$ is
\[
1 \ \ 4 \ \ 10 \ \ 20 \ \ 10 \ \ 4 \ \ 1.
\]
Then Corollary \ref{precise resol} gives the {\em minimal} free
resolution for $R/G$, and (\ref{resol of R/I}) gives the following
$R$-resolution for $R/I$:
\[
0 \rightarrow R(-16)^{25} \rightarrow R(-15)^{52} \rightarrow
\begin{array}{c}
R(-10)^6 \\
\oplus \\
R(-14)^{25}
\end{array}
\rightarrow
\begin{array}{c}
R(-5)^4 \\
\oplus \\
R(-10)
\end{array}
\rightarrow R \rightarrow R/I \rightarrow 0.
\]
So the type of splitting mentioned above does not occur.  The only
possible splitting would be the summand $R(-10)$ at the beginning
of the resolution. However, since we have assumed that $d_{n+1}
\leq (\sum_{i=1}^n d_i) - n$ (Remark \ref{hyp}) and the $G_i$ are
chosen generically, $G_5$ is a minimal generator of $I$ and so
this splitting does not occur either, and the above is the minimal
free $R$-resolution of $R/I$.

In the next section, and in Example \ref{no splitting}, we will see examples
where overlaps arise in other parts of the resolution and still cannot be
split off.
\end{example}

The idea behind Example \ref{ex of no split} leads to one situation where we
can give a minimal free resolution of a general almost complete
intersection in any number of variables.  If one of the generators has degree
1 then we may pass to the analogous problem in a ring of one fewer variable,
so we will assume that our ideal is non-degenerate, i.e.\ that $d_1 \geq 2$.

\begin{corollary}\label{one peak resol}
Let $I = (G_1,\dots,G_{n+1})$ be a general almost complete
intersection in $R = k[x_1,\dots,x_n]$, with $d_i = \deg G_i$ and
$2 \leq d_1 \leq d_2 \leq \dots \leq d_n \leq d_{n+1} \leq
(\sum_{i=1}^n d_i) -n$.  Let $d = d_1+\cdots+d_n$.  Assume that
$(\sum_{i=1}^{n+1} d_i) -n$ is even and that $d_2 + \cdots + d_n <
d_1 + d_{n+1} + n$.  Then $R/I$ has a  free $R$-resolution of the
form
\[
\begin{array}{c}
0 \rightarrow F_1^\vee (-d) \rightarrow
\begin{array}{c}
K_1^\vee (-d) \\ \oplus \\ F_2^\vee (-d)
\end{array}
\rightarrow
\begin{array}{c}
K_2^\vee (-d) \\ \oplus \\ F_3^\vee (-d)
\end{array}
\rightarrow \cdots \rightarrow
\begin{array}{c}
K_{n-2}^\vee (-d) \\ \oplus \\ F_{n-1}^\vee (-d)
\end{array}
\hbox{\hskip 4cm}
\\ \\
\hbox{\hskip 7cm} \rightarrow \bigoplus_{i=1}^{n+1} R(-d_i) \
\rightarrow R \rightarrow R/I \rightarrow 0
\end{array}
\]
where $K_i$ is the $i$-th free module in the Koszul resolution of
$R/(G_1,\dots,G_n)$,
\[
F_i = R(-\ell-i)^{\alpha_i + \gamma_i}  \ \ \hbox{ with }
\
\
\alpha_i = \binom{n+\ell-1}{i+\ell} \binom{\ell-1+i}{\ell} =
\gamma_{n-i},
\]
and $\ell$ is as defined in Remark \ref{observations about MN3}.  If $\ell+1 =
d_i$ for any $1 \leq i \leq n$ then for each such occurrence there is a
corresponding splitting of a free summand from $F_1^\vee (-d)$ and $K_1^\vee
(-d)$.  Except for such splitting, this resolution is minimal.
\end{corollary}

\begin{proof}
As we noted in Remark \ref{observations about MN3}, Lemma \ref{G hilb function
facts} shows that the hypotheses on the $d_i$ are precisely what we need in
order to apply Corollary \ref{precise resol}, and we obtain $s=2\ell$ and
generators of degrees $d_i$, $1 \leq i \leq n+1$.  The resolution was
already observed in (\ref{resol of R/I}).

The splitting of summands of $F_1^\vee (-d)$ when $\ell+1 = d_i$ occurs since
$\ell+1$ is the initial degree of $G$ and $J \subset G$, so any such
polynomial can be taken to be a minimal generator of $G$.  Hence it splits
thanks to a usual mapping cone argument.

Now we show that other overlapping terms do not split.  Because of the way
that the homomorphisms in the mapping cone are formed, there is no splitting
between summands of $K_i^\vee (-d)$ and $F_{i+2}^\vee (-d)$, since there is
no map between these modules.  (See Example \ref{no splitting}.)    There is
no splitting of summands of $K_i^\vee (-d)$ and $K_{i+1}^\vee (-d)$ since the
$K_i$ come from a minimal free resolution, and similarly there is no
splitting of summands of $F_i^\vee (-d)$ and $F_{i+1}^\vee (-d)$.  We now
check that there is no possible overlap between summands of $F_i^\vee (-d)$
and $K_i^\vee (-d)$ for $i \geq 2$ (we have already accounted for overlaps
when $i=1$).

Clearly it is enough to compare summands of $F_i$ and $K_i$.  Each summand of
$F_i$ is of the form $R(-\ell-i)$, while the summands of $K_i$ are of the form
$R(-d_{r_1} - \cdots - d_{r_i})$.  So we have to show that it is impossible
to have an equality of the form
\[
\frac{d-d_{n+1} - n}{2} + i = d_{r_1} + \cdots + d_{r_i}
\]
for $i \geq 2$, i.e.\ it is impossible to have $d-d_{n+1} - n + 2i = 2d_{r_1}
+ \cdots + 2d_{r_i}$.  Let $A > 0$ be the integer such that $d_2 + \cdots +
d_n = d_1 + d_{n+1} + n - A$.  Then we have
\[
\begin{array}{rcl}
d-d_{n+1} -n+2i & = & (d_1 + d_{n+1} + n - A) + d_1 - d_{n+1} - n + 2i \\
& = & 2d_1 + 2i - A \\
& = & 2d_{r_1} + \cdots + 2d_{r_i} + (2d_1 + 2i - A - 2d_{r_1} - \cdots -
2d_{r_i}).
\end{array}
\]
So the suggested equality holds if and only if
\[
\begin{array}{rcl}
A & = & 2d_1 + 2i -2d_{r_1} - \cdots - 2d_{r_i} \\
& = & 2[(d_1 - d_{r_1}) + (i - d_{r_2} - \cdots - d_{r_i})]
\end{array}
\]
If $i=1$ this could be strictly positive (if and only if $d_1 = d_{r_1}$), as
we have seen.  If
$i=2$ this could be non-negative (if and only if $i=2$ and $d_1 = d_{r_1} =
d_{r_2} = 2$) but not strictly positive.  If $i \geq 3$ this is strictly
negative, and we have our contradiction.
\end{proof}

\begin{example} \label{no splitting}
Note that there can be ``ghost'' terms that do not split off,
coming from overlaps between $K_i^\vee (-d)$ and $F_{i+2}^\vee
(-d)$.  For instance, let $n=4$, $d_1 = 3$, $d_2 = d_3 =d_4 = 5$,
$d_5 = 10$.  Then one computes $\ell=2$ and $d=18$, and
(\ref{resol of R/I}) predicts a  free $R$-resolution
\[
0 \rightarrow R(-15)^{16} \rightarrow
\begin{array}{c}
\left [
\begin{array}{c}
R(-15)  \\
\oplus \\
R(-13)^3
\end{array}
\right ]
\\
\oplus \\
R(-14)^{30}
\end{array}
\rightarrow
\begin{array}{c}
\left [
\begin{array}{c}
R(-10)^3 \\
\oplus \\
R(-8)^3
\end{array}
\right ] \\
\oplus \\
R(-13)^{16}
\end{array}
\rightarrow
\begin{array}{c}
\left [
\begin{array}{c}
R(-3) \\
\oplus \\
 R(-5)^3
\end{array}
\right ] \\
\oplus \\
R(-10)
\end{array}
\rightarrow R \rightarrow R/I \rightarrow 0.
\]
Our observation above allows a splitting off of a summand $R(-15)$, and
Corollary \ref{one peak resol} guarantees that the rest of the resolution is
minimal, despite the overlap of three copies of
$R(-13)$ in the second and third free modules in the resolution and the
overlap of one copy of $R(-10)$ in the first and second.  (We have already
seen in Example \ref{ex of no split} that the summand $R(-10)$ in the first
module does not split off.)
\end{example}

The hypotheses of Corollary \ref{one peak resol} guarantee that
the Hilbert function of $R/G$ has one peak and grows maximally,
respectively (Lemma \ref{G hilb function facts} (b) and (d)).
These are enough to completely determine the minimal free
$R$-resolution.  We would like to weaken these hypotheses as much
as possible.  One way is Corollary \ref{weak bound}, where these
hypotheses are removed completely but a bound for the Betti
numbers of $R/G$ is used which greatly exceeds those which appear
in our ``general'' situation.

We will make a slight improvement, allowing two peaks (by Lemma
\ref{G hilb function facts} no more peaks can occur) but keeping
the fact that the Hilbert function of $R/G$ grows maximally . Note
first that in this situation the minimal free resolution is not
uniquely determined just from the Hilbert function.

\begin{example}
Consider the case $n=4$, $d_1=d_2=d_3=d_4=4$, $d_5=5$.  The corresponding
Gorenstein Hilbert function is
\[
1 \ \ 4 \ \ 10 \ \ 20 \ \ 20 \ \ 10 \ \ 4 \ \ 1.
\]
 Thanks to \cite{MN3} the maximal possible graded Betti numbers are
\[
0 \rightarrow R(-11) \rightarrow
\begin{array}{c}
R(-6)^{10} \\
\oplus \\
R(-7)^{15}
\end{array}
\rightarrow
\begin{array}{c}
R(-5)^{24} \\
\oplus \\
R(-6)^{24}
\end{array}
\rightarrow
\begin{array}{c}
R(-4)^{15} \\
\oplus \\
R(-5)^{10}
\end{array}
\rightarrow R \rightarrow R/G \rightarrow 0,
\]
while thanks to \cite{boij} the general such graded Betti numbers are
\[
0 \rightarrow R(-11) \rightarrow
R(-7)^{15}
\rightarrow
\begin{array}{c}
R(-5)^{14} \\
\oplus \\
R(-6)^{14}
\end{array}
\rightarrow
R(-4)^{15}
\rightarrow R \rightarrow R/G \rightarrow 0.
\]
\end{example}

So, we first need to determine the generic Betti numbers of certain 
Gorenstein Artinian graded algebras. In \cite{boij}, Corollary 3.10 Boij
established the existence of generic Betti numbers of compressed level
algebras and in Conjecture 3.13 he guessed that the generic Betti numbers
are as small as we can hope for.  Using the fact that the Minimal Resolution
Conjecture (MRC) has been proved for large numbers of points in any $\PP^r$
by Hirschowitz and Simpson (\cite{HS}) and that the canonical module of the
homogeneous coordinate ring of points can be identified with an
ideal of the ring itself, we will  prove (resp.\ partially
prove) that Boij's Conjecture (\cite{boij}, Conjecture 3.13)
holds for Gorenstein Artinian algebras of even
(resp.\ odd) embedded dimension $n$, initial degree $t\gg 0$
and socle degree $2t-1$.

\begin{proposition}\label{gor two peaks even}  Let $A$ be a generic
 Gorenstein Artinian graded algebra  of even embedding dimension 
$n=2p$ with initial degree $t$ and socle degree $2t-1$. Assume
 $n=4$, or $n>4$ and $t\gg 0$.  Then $A$ has a minimal free
 $R$-resolution of the following type:
\[
0 \rightarrow R(-2t-n+1)\rightarrow
R(-t-n+1)^{\alpha_1} \rightarrow ... \rightarrow R(-t-p-1)^{\alpha
_{p-1}} \rightarrow \hbox{\hskip 1in}
\]
\[
 \left (
\begin{array}{c} \displaystyle
 R(-t-p)^{\alpha_{p}} \\ \oplus \\ R(-t-p+1)^{\alpha_{p}}
\end{array}
\right ) \rightarrow
 R(-t-p+2)^{\alpha _{p-1}}  \rightarrow .... \rightarrow
  R(-t-1)^{\alpha_2} \rightarrow 
\]
\[
 \hskip 2in R(-t)^{\alpha_1} \rightarrow  R 
\rightarrow A \rightarrow 0
\]
where 
\[
\alpha_{i}=\binom{t+n-1}{t+i-1}\binom{t+i-2}{i-1}
-\binom{t+n-1}{ t+n-i}\binom{t+n-i-1}{n-i}
\]
 for $i=1,...p$.

\end{proposition}

\begin{proof} 
For $n=4$ see \cite{boij}, Proposition 3.24. Assume
$n=2p>4$ and $t\gg 0$.
 According to Hirschowitz and Simpson (\cite{HS}; Th\'eor\`eme) there 
exists an integer $\pi (n-1)$ ($\pi (n-1) \sim 6^{(n-1)^3log(n-1)}$)
such that the MRC holds for $s> \pi(n-1)$ points in general position in
$\PP^{n-1}$. Let $X\subset \PP^{n-1}$ be a set of 
\[
\rho
(t,n):=\left \lceil
\frac{t(t+1)...(t+p-1)(t+p+1)...(t+n-1)}{(n-1)!}\right \rceil 
\]
points in generic position. (For any $x\in \RR$, we set $\lceil
x\rceil := min\{n\in \ZZ \mid  x\le n\}$). Since $t\gg 0$, we have
$\rho (t,n)> \pi(n-1)$.

\vskip 2mm We first observe that $$\binom{ t+n-2}{n-1}\le \rho
(t,n)=\left \lceil\frac
{t(t+1)...(t+p-1)(t+p+1)...(t+n-1)}{(n-1)!}\right \rceil<\binom{
t+n-1}{n-1}.$$

\vskip 2mm \noindent In particular, $R/I(X)$ has generic Hilbert
function with initial degree $t$. As a result, $I(X)$ has a
minimal free $R=k[x_0,...,x_{n-1}]$-resolution of the following
type:
\[
\begin{array}{c}
0 \rightarrow \left (
\begin{array}{c} \displaystyle
 R(-t-n+1)^{b_{n-1}} \\ \oplus \\ R(-t-n+2)^{a_{n-1}}
\end{array}
\right ) \rightarrow \left (
\begin{array}{c} \displaystyle
 R(-t-n+2)^{b_{n-2}} \\ \oplus \\ R(-t-n+3)^{a_{n-2}}
\end{array}
\right ) \rightarrow .... \rightarrow

 \hbox{\hskip 2cm}
\\ \\
\hbox{\hskip 3cm} \displaystyle
  \left (
\begin{array}{c} \displaystyle
 R(-t-2)^{b_{2}} \\ \oplus \\ R(-t-1)^{a_{2}}
\end{array}
\right ) \rightarrow
 \left (
\begin{array}{c} \displaystyle
 R(-t-1)^{b_{1}} \\ \oplus \\ R(-t)^{a_{1}}
\end{array}
\right ) \rightarrow R \rightarrow R/I(X) \rightarrow 0

\end{array}
\]
where 
\[
a_i = \max \{ 0,
h^0(\PP^{n-1},\Omega^{i-1}_{\PP^{n-1}}(t+i-1))- rk(\Omega
^{i-1}_{\PP^{n-1}})\rho(t,n) \}
\]
 and 
\[
b_i = \max \{ 0, rk(\Omega^{i}_{\PP^{n-1}})
\rho(t,n)-h^0(\PP^{n-1},\Omega^{i}_{\PP^{n-1}}(t+i)) \}
\]
 for any $i=1,...,n-1$ (See \cite{HS}; pg. 2). Note that we set
$ \Omega^{0}_{\PP^{n-1}}={\cO}_{\PP^{n-1}}$.

\vskip 2mm 

An intricate calculation shows that for any $i$ we have
\[
\begin{array}{rcl}
a_i & = & \max\{0, h^0(\PP^{n-1},\Omega^{i-1}_{\PP^{n-1}}(t+i-1))-
rk(\Omega ^{i-1}_{\PP^{n-1}})\rho(t,n)\} \\ \\

& = &\displaystyle  \max \left \{ 0,\binom{t+n-1
}{ t+i-1}\binom{t+i-2 }{i-1} -\binom{n-1 }{i-1}\rho (t,n) \right \}
\end{array}
\]
and
\[
\begin{array}{rcl}
b_i & = & \max \{0,rk(\Omega
^{i}_{\PP^{n-1}})\rho(t,n)
- h^0(\PP^{n-1},\Omega^{i}_{\PP^{n-1}}(t+i)) \} \\ \\
& = & \displaystyle \max \left \{ 0, \binom{n-1 }{ i}\rho
(t,n)-\binom{t+n-1}{ t+i} \binom{t+i-1}{ i} \right \}
\end{array}
\]
Furthermore,

\begin{itemize}
\item if  $p+1\le i \le n-1$ then $a_i = 0$
\item  if $1\le i \le p-1$ then $b_i = 0$.
\end{itemize}

\vskip 2mm  Hence, the minimal free $R$-resolution of $R/I(X)$ has
the form

\[
0 \rightarrow
 R(-t-n+1)^{b_{n-1}}  \rightarrow ..... \rightarrow
 R(-t-p-1)^{b_{p+1}}
\rightarrow \begin{array}{c}R(-t-p)^{b_{p}}\\ \oplus \\
R(-t-p+1)^{a_{p}} \end{array}
\]
\[
 \rightarrow R(-t-p+2)^{a_{p-1}}  \rightarrow ..... \rightarrow
  R(-t-1)^{a_{2}} \rightarrow  R(-t)^{a_{1}}
 \rightarrow R \rightarrow R/I(X) \rightarrow 0
\]
where $$a_i=\binom{t+n-1 }{ t+i-1}\binom{t+i-2 }{i-1} -\binom{n-1
}{i-1}\rho (t,n)$$and $$b_i=\binom{n-1 }{ i}\rho
(t,n)-\binom{t+n-1}{ t+i}\binom{t+i-1}{ i}. $$

Let $A(X)=R/I(X)$ be the homogeneous coordinate ring of $X$. The
canonical module $\omega_X$ of $A(X)$ can be embedded as an ideal
$\omega_X \subset A(X)$ of initial degree $t$ and we have a short
exact sequence 
\[
0 \rightarrow \omega _X \rightarrow A(X)
\rightarrow A \rightarrow 0
\]
 where $A$ is a Gorenstein Artinian algebra of embedding dimension $n$
with initial degree $t$ and socle degree $2t-1$ (cf.\ \cite{boij2}). A
straightforward calculation shows that 
\[
\begin{array}{rcl}
a_{i} + b_{n-i} & =
& h^0(\PP^{n-1},\Omega^{i-1}_{\PP^{n-1}}(t+i-1))
-h^0(\PP^{n-1},\Omega^{i}_{\PP^{n-1}}(t+i)) \\ \\

& = & \displaystyle \binom{t+n-1 }{ t+i-1}\binom{t+i-2}{ i-1}
-\binom{t+n-1}{ t+n-i}\binom{t+n-i-1 }{n-i}.
\end{array}
\]
 Hence, applying the mapping cone
 construction to the diagram
\[
\begin{array}{ccccccccccccccc} &&&&&&&&&& 0 \\
&&&&&&&&&& \downarrow
\\
0 & \rightarrow & R(-2t-n+1) & \rightarrow & R(-t-n+1)^{a_1} &
\rightarrow & ... & \rightarrow
 & R(-t)^{b_{n-1}} & \rightarrow
 & \omega _X &
\rightarrow & 0 \\ && \downarrow && \downarrow
 &&      && \downarrow  && \downarrow \\
0  & \rightarrow  & R(-t-n+1)^{b_{n-1}} & \rightarrow &
R(-t-n+2)^{b_{n-2}}& \rightarrow &
 ... &  \rightarrow & R &
 \rightarrow & R/I(X) & \rightarrow & 0 \\&
&&&& && &&&  \downarrow \\ &&&&&&&&&& A \\ &&&&&&&&&& \downarrow
\\
 &&&&&&&&&& 0

\end{array}
\]
 we get the
minimal free $R$-resolution of $A$:

\[0 \rightarrow R(-2t-n+1)\rightarrow
R(-t-n+1)^{\alpha_1} \rightarrow ... \rightarrow R(-t-p-1)^{\alpha
_{p-1}} \rightarrow
\]
\[
\begin{array}{c}
 \left (
\begin{array}{c} \displaystyle
 R(-t-p)^{\alpha_{p}} \\ \oplus \\ R(-t-p+1)^{\alpha_{p}}
\end{array}
\right ) \rightarrow
 R(-t-p+2)^{\alpha _{p-1}}  \rightarrow .... \rightarrow
  R(-t-1)^{\alpha_2} \rightarrow R(-t)^{\alpha_1} \rightarrow R \rightarrow
A \rightarrow 0
\end{array}
\]
with 
\[
\alpha_{i} = \binom{t+n-1}{t+i-1}\binom{t+i-2}{i-1}
-\binom{t+n-1}{ t+n-i}\binom{t+n-i-1}{n-i}
\]
 for $i=1,...p$.
\end{proof}

\begin{remark}\label{gor two peaks odd} Arguing as above we can prove that the
 generic minimal free $R$-resolution  of a Gorenstein Artinian
 graded
algebra $A$ of odd embedding dimension $n = 2p+1>3$ with initial
degree $t \gg 0$ and socle degree $2t-1$ has the form
\[
\begin{array}{c}
0 \rightarrow R(-2t-n+1)\rightarrow R(-t-n+1)^{\alpha _1} \rightarrow
... \rightarrow R(-t-p-1)^{\alpha _p}\oplus R(-t-p)^{\le b_p}
\\ \\
\rightarrow R(-t-p+1)^{\alpha _p}\oplus R(-t-p)^{\le b_p}\rightarrow
\dots \rightarrow  R(-t-1)^{\alpha _2} \rightarrow R(-t)^{\alpha _1}
\rightarrow R \rightarrow A \rightarrow 0
\end{array}
\]
with 
\[
\alpha _{i} = \binom{t+n-1}{t+i-1}\binom{t+i-2}{i-1}
-\binom{t+n-1}{ t+n-i}\binom{t+n-i-1}{ n-i}
\]
 for $i=1,...,p$ and
\[
b_p = \binom{n-1 }{ p}\rho (t,n)-\binom{t+n-1}{ t+p}\binom{t+p-1}{
p} 
\]
 but we do not know if the overlap appearing in the middle
of the resolution does occur.
\end{remark}

We have already found the minimal free resolution when there is
one peak and ``maximal" Hilbert function (Corollary \ref{one peak
resol}), so without loss of generality we will assume two peaks (i.e.\
$\sum_{i=1}^{n+1} d_i -n$ is odd).  The case $n=3$ will be treated
in the next section.  Now we give the analogous result when $n$ is even
and the socle degree is large.  We have

\begin{proposition}\label{aci two peaks}
Let $I = (G_1,\dots,G_{n+1})$ be a general almost complete
intersection in $R = k[x_1,\dots,x_n]$, with $d_i = \deg G_i$, $2
\leq d_1 \leq d_2 \dots \leq d_n \le d_{n+1} \leq (\sum_{i=1}^n d_i)
-n$, $\sum_{i=1}^{n+1} d_i -n$  odd and $n=2p$. Let
\[
\begin{array}{rcl}
d & = & d_1+ \dots +d_n \\ \\ \ell &= & \displaystyle \frac{d-d_{n+1}
-n-1}{2}.
\end{array}
\]
Assume that $n=4$ and $d_2+d_3+d_4<d_1+d_5+4$ or $n>4$,  $d_2 +
 \dots  +d_ n < d_1 + d_{n+1} + n$ and $\ell \gg 0$.
Then $R/I$ has a free $R$-resolution of the form

\[
\begin{array}{c}
0 \rightarrow F_1^\vee (-d) \rightarrow
\begin{array}{c}
K_1^\vee (-d) \\ \oplus \\ F_2^\vee (-d)
\end{array}
\rightarrow
\begin{array}{c}
K_2^\vee (-d) \\ \oplus \\ F_3^\vee (-d)
\end{array}
\rightarrow \cdots \rightarrow
\begin{array}{c}
K_{n-2}^\vee (-d) \\ \oplus \\ F_{n-1}^\vee (-d)
\end{array}
\hbox{\hskip 4cm}
\\ \\
\displaystyle
\hbox{\hskip 7cm} \rightarrow \bigoplus_{i=1}^{n+1} R(-d_i) \
\rightarrow R \rightarrow R/I \rightarrow 0
\end{array}
\]
where $K_i$ is the $i$-th free module in the Koszul resolution of
$R/(G_1,\dots,G_n)$ and

\[ F_i = \left \{
\begin{array}{ll}
R(-\ell-i)^{\alpha_i} & \hbox{if $1 \leq i \leq p-1$};
\\R(-\ell-p)^{\alpha_p} \oplus R(-\ell-p-1)^{\alpha_p}& \hbox{if $i=p$};
\\R(-\ell-i-1)^{\alpha_{n-i}}  & \hbox{if $p+1 \leq i \leq n-1$}
\end{array}
\right.
\]
where 
\[
\alpha_i =\binom{\ell +n }{ \ell +i}\binom{\ell +i-1 }{i-1}
-\binom{\ell +n }{\ell +1+n-i}\binom{\ell +n-i }{n-i}
\]
 for $i=1,...,p$.

\vskip 2mm

\begin{itemize}
\item[a.]   If $\ell+1 = d_i$ for any $1 \leq i \leq n$ then for each such
occurrence there is a corresponding splitting of a free summand at
the end of the resolution.

\item[b.] If $n=4$, $d_1=d_2=2$ and $d_3+d_4 = d_5+3$
then there is a splitting of one term of the form $R(\ell +3-d)$
with $R(-d_3-d_4)$.
\end{itemize}

Except for such splitting, this resolution is minimal.

\end{proposition}
\begin{proof}
The calculations are almost identical with those of Corollary
\ref{one peak resol}, but there is a difference in the reasoning.
Here we will prove the existence of an almost complete
intersection with the claimed resolution, and we will find all
possible splittings that can occur numerically, and show that
these do occur.  After such splitting, by semicontinuity, the
result must be the minimal free $R$-resolution of a general almost
complete intersection.

 By Proposition \ref{gor two peaks even} for a {\em general}
Gorenstein Artinian quotient of $k[x_1, ... ,x_{n=2p}]$ whose
Hilbert function is maximal (i.e.\ compressed) and with two peaks,
the minimal free $R$-resolution is of the following form:
\[ 
\begin{array}{c}
0 \rightarrow R(-2\ell -n-1)\rightarrow
R(-\ell-n)^{\alpha _1} \rightarrow ... \rightarrow R(-\ell
-p-2)^{\alpha _{p-1}}
 \rightarrow \hbox{\hskip 1in} 
\\ \\
 R(-\ell -p-1)^{\alpha _p}  \oplus
 R(-\ell -p)^{\alpha _p}
 \rightarrow
 R(-\ell -p+1)^{\alpha
_{p-1}}\rightarrow... \rightarrow  R(-\ell -2)^{\alpha _2} 
\\ \\
\hbox{\hskip 3in} \rightarrow R(-\ell -1)^{\alpha _1} \rightarrow R
\rightarrow R/G
\rightarrow 0
\end{array}
\]
where 
\[
\alpha _{i} = \binom{\ell +n }{ \ell +i}\binom{\ell +i-1
}{i-1} -\binom{\ell +n }{ \ell + n+1-i}\binom{\ell +n-i}{ n-i}
\]
 for $i=1,...,p$. (Note that there is a slight difference of notation:
the value of $t$ in Proposition \ref{gor two peaks even},  is now
$\ell +1$.)

Now, notice that such a general Gorenstein Artinian quotient has
all its generators in degree $\ell +1$, and a simple calculation
(using the hypotheses in the statement of the Proposition) shows
that $\ell +1 \leq d_1$.  Hence a complete intersection $J \subset
G$ exists with generators of degrees $d_1,...,d_n$.  By a standard
mapping cone argument, the residual $I=[J:G]$ is an almost
complete intersection (Cohen-Macaulay of height $n$ with $n+1$
minimal generators).  The discussion above shows that $I$ has a
free resolution of the form claimed.

Now we consider splitting. The type of splitting mentioned in part
a.\ has already been discussed.  As we have already observed, the
only possible splitting comes (in the resolution (\ref{resol of
R/I})) between a summand of $F_i^\vee (-d)$ and one of $K_i^\vee
(-d)$. So  let us study the possible overlaps between summands of
$F_i^\vee (-d)$ and of $K_i^\vee (-d)$ for $i\ge 2$.  (We have
already accounted in part a.\ for  overlaps when $i=1$.)  First we
make a numerical calculation.  \medskip

\underline{\em Claim:} {\em Under the hypotheses of part b.\ there
is exactly one summand in common between $K_2^\vee (-d)$ and
$F_2^\vee (-d)$.  If these hypotheses are not met then there is no
other summand in common.}

Arguing as in Corollary \ref{one peak resol} we check that if
$n=2p$ and $2\le i \le p$ there is no overlapping between summands
of $R(\ell +i-d)^{\alpha _i} $ and summands of 
\[
K_i^\vee
(-d) = \bigoplus_{d_1\le d_{r_1}<...<d_{r_i}\le
d_n}R(d_{r_1}+...+d_{r_i}-d).
\]
 Assume now that  $n=2p$ and $p\le i
\le n-1$. There is a summand in common between $R(\ell
+i+1-d)^{\alpha _{n-i}} $ and $K_i^\vee (-d)=\oplus _{d_1\le
d_{r_1}<...<d_{r_i}\le d_n}R(d_{r_1}+...+d_{r_i}-d)$ if and only
if $\ell + i +1= d_{r_1} + \cdots + d_{r_i}$.
 Let $A > 0$ be the integer such
that $d_2 + \cdots + d_n = d_1 + d_{n+1} + n - A$.  We have 
\[
\begin{array}{rcl}
 A  & = & 2d_1 + 2i -2d_{r_1} - \cdots - 2d_{r_i} \\ 
& = & 2(d_1 - d_{r_1}) + 2(i - d_{r_2} - \cdots - d_{r_i})+1. 
\end{array}
\]
 So, necessarily, we have $i=2$ ( and hence $n=4$),
$d_{r_1}=d_{r_2}=2$. Since we have assumed $d_1>1$ this forces
$n=4$, $d_1=d_2=2$ and $d_3+d_4=d_5+3$ as desired.

So after the splitting of part a.\ we have a minimal free
resolution {\em unless} the numerical conditions of b.\ hold.  In
this case we have to prove one more splitting.

Assume that the conditions of b.\ hold.  Then one checks that
$\ell = 1$, so the Hilbert function of $R/G$ is $(1 \ \ 4 \ \ 4 \
\ 1)$ and Proposition \ref{gor two peaks even} gives that $G$ is
generated in degree 2 and has 5 syzygies in degree 4.  In fact, the minimal
free resolution of $R/G$ is
\begin{equation}\label{desired resol}
0 \rightarrow R(-7) \rightarrow R(-5)^6 \rightarrow R(-4)^5 \oplus R(-3)^5
\rightarrow R(-2)^6 \rightarrow R \rightarrow R/G \rightarrow 0.
\end{equation}
One can check that the generators of $I$ are of degrees 2, 2, $d_3$, $d_4$,
$d_5$ where $d_3 >2$.  Then $J$ also has exactly two generators of degree 2,
and hence exactly one first syzygy of degree 4.  It is the corresponding
summand of $K_2^\vee (-d)$ that we would like to split with a summand of
$F_2^\vee (-d)$.  As before, our strategy will be to construct a specific
$G$ and link which meets our needs, and then the same will hold for the
generic case.  

Let $Z \subset \proj{3}$ be a general set of 5 points.  Then $Z$ is
arithmetically Gorenstein with Hilbert function
\[
1 \ \ 4 \ \ 5 \ \ 5 \ \ \rightarrow
\]
and minimal free resolution
\[
0 \rightarrow R(-5) \rightarrow R(-3)^5 \rightarrow R(-2)^5 \rightarrow R
\rightarrow R/I_Z \rightarrow 0.
\]
 Let $Q$ be a generally chosen form of degree 2.  Then $G := I_Z + (Q)$ is
the saturated ideal of a height 4 Gorenstein ideal with Hilbert function
$(1 \ \ 4 \ \ 4 \ \ 1)$ and minimal free resolution (\ref{desired resol})
(obtained by a tensor product).  Note that $G$ has five linear syzygies and
five quadratic syzygies, and that $I_Z$ has five linear syzygies.  Hence the
five linear syzygies of $G$ are precisely the five linear syzygies of $I_Z$,
and the five quadratic syzygies of $G$ are the Koszul syzygies of $Q$ with
the five generators of $I_Z$.  Now choose $J$ to consist of $Q$, one quadric 
generator of $I_Z$, and then general forms in $G$ of suitable degrees
satisfying the numerical conditions (e.g.\ $d_3 = d_4 = 4$, $d_5 = 5$). 
This produces an $I$ as desired, and clearly the Koszul syzygy of degree 4
for $J$ is the one from $G$.

Since no further numerical overlap exists, the remaining
resolution is minimal, and hence is the minimal free
$R$-resolution of the general almost complete intersection of type
$(d_1,...,d_{n+1})$.
\end{proof}

\begin{example}  \label{codim 4 result} Suppose $n=4$.
 Let $I = (G_1,\dots,G_{5})$ be a general
almost complete intersection in $R = k[x_1,\dots,x_4]$, with $d_i
= \deg G_i$, $2 \leq d_1 \leq d_2 \leq d_3 \leq d_4 \leq d_5 \leq
(\sum_{i=1}^4 d_i) -4$ and $\sum_{i=1}^5 d_i -4$  odd. Let
$d  =  d_1+d_2+d_3+d_4$ and $\ell =  \frac{d -d_5 -5}{2}$.
Assume that  $d_2 + d_3 + d_4 < d_1 + d_5 + 4$.
 Then $R/I$ has a  free $R$-resolution of the form
\[
\begin{array}{c}
0 \rightarrow R(\ell+1-d)^a \rightarrow \left (
\begin{array}{c} \displaystyle
\bigoplus_{i=1}^4 R(d_i-d) \\ \oplus \\ R(\ell +2-d)^{b} \\ \oplus
\\ R(\ell+3-d)^{b}
\end{array}
\right ) \rightarrow \left (
\begin{array}{c} \displaystyle
\bigoplus_{1\leq i<j \leq 4} R(-d_i-d_j) \\ \oplus \\
R(\ell+4-d)^a
\end{array}
\right ) \hbox{\hskip 2cm}
\\ \\
\hbox{\hskip 7cm} \displaystyle \rightarrow \bigoplus_{i=1}^5
R(-d_i) \rightarrow R \rightarrow R/I \rightarrow 0
\end{array}
\]
where we have $a = \binom{\ell+3}{2}$ and $b  =  \binom{\ell+3}{2}
-1 $.
\begin{itemize}
\item[a.]   If $\ell+1 = d_i$ for any $1 \leq i \leq 4$ then for each such
occurrence there is a corresponding splitting of a free summand at
the end of the resolution.

\item[b.] If $d_1=d_2=2$ and $d_3+d_4 = d_5+3$  then there is a splitting of
one term of the form $R(\ell +3-d)$ with $R(-d_3-d_4)$.
\end{itemize}

Except for such splitting, this resolution is minimal.

\end{example}

\begin{remark}\label{aci two peaks n odd}
Let $I = (G_1,\dots,G_{n+1})$ be a general almost complete
intersection in $R = k[x_1,\dots,x_n]$, with $d_i = \deg G_i$, $2
\leq d_1 \leq d_2 ... \leq d_n \le d_{n+1} \leq (\sum_{i=1}^n d_i)
-n$, $\sum_{i=1}^{n+1} d_i -n$  odd and $n=2p+1>3$. Let
\[
\begin{array}{rcl}
d & = & d_1+...+d_n \\ \\ \ell &= & \displaystyle \frac{d-d_{n+1}
-n-1}{2}.
\end{array}
\]
Assume that   $d_2 + ... +d_ n < d_1 + d_{n+1} + n$ and $\ell \gg
0$.  Then $R/I$ has a free $R$-resolution of the form

\[
\begin{array}{c}
0 \rightarrow F_1^\vee (-d) \rightarrow
\begin{array}{c}
K_1^\vee (-d) \\ \oplus \\ F_2^\vee (-d)
\end{array}
\rightarrow
\begin{array}{c}
K_2^\vee (-d) \\ \oplus \\ F_3^\vee (-d)
\end{array}
\rightarrow \cdots \rightarrow
\begin{array}{c}
K_{n-2}^\vee (-d) \\ \oplus \\ F_{n-1}^\vee (-d)
\end{array}
\hbox{\hskip 4cm}
\\ \\
\hbox{\hskip 7cm} \rightarrow \bigoplus_{i=1}^{n+1} R(-d_i) \
\rightarrow R \rightarrow R/I \rightarrow 0
\end{array}
\]
where $K_i$ is the $i$-th free module in the Koszul resolution of
$R/(G_1,\dots,G_n)$ and
\[ 
F_i = \left \{
\begin{array}{ll}
R(-\ell-i)^{\alpha_i} & \hbox{if $1 \leq i \leq p-1$};
\\R(-\ell-p)^{\alpha_p} \oplus R(-\ell-p-1)^{\le b_p}& \hbox{if $i=p$};
\\R(-\ell-p-2)^{\alpha_p} \oplus R(-\ell-p-1)^{\le b_p}& \hbox{if $i=p+1$};
\\R(-\ell-i-1)^{\alpha_{n-i}}  & \hbox{if $p+2 \leq i \leq n-1$}
\end{array}
\right.
\]
where 
\[
\alpha_i =\binom{\ell +n }{ \ell +i}\binom{\ell +i-1 }{i-1}
-\binom{\ell +n }{\ell +1+n-i}\binom{\ell +n-i }{n-i}
\]
 for $i=1,...,p$ and 
\[
b_p = \binom{n-1 }{ p}\rho (t,n)-\binom{t+n-1}{
t+p}\binom{t+p-1}{ p}. 
\]
 Moreover, if $\ell+1 = d_i$ for any $1
\leq i \leq n$ then for each such occurrence there is a
corresponding splitting of a free summand at the end of the
resolution.
\end{remark}

%%%%%%%%%%%%%%%%%%%%%%%%%%%%%%%%%%%%%%%%%%%%%%%%%%%%%%%%%%%%%%%%%%%%%%%%%%%%%

\section{The case $n=3$}

We now specialize to the case $n=3$, and for simplicity of notation we will
write $R = k[x,y,z]$.  The goal of this section is to find the minimal free
resolution for any general almost complete intersection $I = (G_1,G_2,G_3,G_4)$
in $R$, with generators of arbitrary degrees $d_i = \deg G_i$.  As before we
assume that $d_1 \leq d_2 \leq d_3 \leq d_4 \leq d_1+d_2+d_3-3$ (Remark
\ref{hyp}).

We begin with
the resolution (\ref{resol of R/I}), which in our context now becomes
\begin{equation} \label{codim 3 resol}
0 \rightarrow F_1^\vee (-d) \rightarrow
\left (
\begin{array}{c}
\displaystyle
\bigoplus_{i=1}^3 R(d_i -d) \\
\oplus \\
F_2^\vee (-d)
\end{array}
\right )
\rightarrow \bigoplus_{i=1}^4 R(-d_i) \rightarrow R \rightarrow R/I
\rightarrow 0
\end{equation}
where  $d = d_1+d_2+d_3$, and the $F_i$ come from the minimal free
$R$-resolution of $R/G$:
\[
0 \rightarrow R(-e) \rightarrow F_2 \rightarrow F_1 \rightarrow R \rightarrow
R/G \rightarrow 0
\]
with $e = d_1 + d_2 + d_3 - d_4$.

Our approach is similar to that of Proposition \ref{aci two
peaks}. We know the Hilbert function of $R/G$, from which we can
calculate all possible minimal free resolutions using
\cite{diesel}.  The general such can be determined, and a link
gives the resolution of the general almost complete intersection.

\begin{proposition} \label{resol of general G}
Let $I = (G_1,G_2,G_3,G_4)$ be a general almost complete intersection in $R$ of
type $(d_1,d_2,d_3,d_4)$ with $d_1 \leq d_2 \leq d_3 \leq d_4 \leq
d_1+d_2+d_3-3$. Let $G$ be the Gorenstein ideal linked to $I$ by the complete
intersection $J = (G_1,G_2,G_3)$.
Let
\[
s = d_1+d_2+d_3-d_4-3 \ \hbox{ and } \ \ell = \left \lfloor \frac{s}{2} \right
\rfloor.
\]

\begin{itemize}
\item[a.] We have maximal Hilbert function
\[
h_{R/G}(t) = \binom{t+2}{2} \hbox{ for all integers } 0 \leq t \leq
\ell
\]
if and only if $d_2+d_3 < d_1+d_4 + 3$.  In any case the Hilbert function is
described in Lemma \ref{G hilb function facts}.

\item[b.] We have $d_2 > \ell +1$.

\item[c.] For the {\em general} Gorenstein ideal $G'$ with the Hilbert function
described in a., the minimal free resolution of $G'$ is described as follows.

\begin{itemize}
\item[Case I.] \underline{$d_2+d_3 < d_1+d_4 + 3$ and $d_1+d_2+d_3+d_4$ odd:}
\[
0 \rightarrow R(-2\ell-3) \rightarrow R(-\ell-2)^{2\ell +3} \rightarrow R(-\ell
-1)^{2\ell +3} \rightarrow R \rightarrow R/G' \rightarrow 0.
\]

\item[Case II.] \underline{$d_2+d_3 < d_1+d_4 + 3$ and $d_1+d_2+d_3+d_4$ even:}
\[
0 \rightarrow R(-2\ell-4) \rightarrow
\begin{array}{c}
R(-\ell-3)^{\ell+2} \\
\oplus \\
R(-\ell-2)^\delta
\end{array}
\rightarrow
\begin{array}{c}
R(-\ell-1)^{\ell+2} \\
\oplus \\
R(-\ell-2)^\delta
\end{array}
\rightarrow R \rightarrow R/G' \rightarrow 0
\]
where $\delta = 1$ if $\ell$ is even and 0 otherwise.

\medskip

\item[Case III.] \underline{$d_2+d_3 \geq d_1+d_4 + 3$ and $d_1+d_2+d_3+d_4$
odd:}
\[
0 \rightarrow R(-2\ell-3) \rightarrow
\begin{array}{c}
R(-2\ell-3+d_1) \\
\oplus \\
R(-\ell-2)^{2d_1}
\end{array}
 \rightarrow
\begin{array}{c}
R(-d_1) \\
\oplus \\
R(-\ell-1)^{2d_1}
\end{array}
\rightarrow R \rightarrow R/G' \rightarrow 0.
\]

\item[Case IV.] \underline{$d_2+d_3 \geq d_1+d_4 + 3$ and $d_1+d_2+d_3+d_4$
even:}
\[
0 \rightarrow R(-2\ell-4) \rightarrow \begin{array}{c}
R(-2\ell-4+d_1) \\
\oplus \\
R(-\ell -3)^{d_1} \\
\oplus \\
R(-\ell-2)^\delta
\end{array}
\rightarrow
\begin{array}{c}
R(-d_1) \\
\oplus \\
R(-\ell-1)^{d_1} \\
\oplus \\
R(-\ell-2)^\delta
\end{array}
\rightarrow R \rightarrow R/G' \rightarrow 0
\]
where $\delta = 1$ if $d_1$ is odd and 0 otherwise.

\end{itemize}

\end{itemize}

\end{proposition}

\begin{proof}
Note that by Lemma  \ref{G hilb function facts}, the knowledge of the Hilbert
function up to degree $\ell$ determines the full Hilbert function of
$R/G$ by symmetry.  The same lemma then gives the values of $h_{R/G}(t)$,
completing part a.  Part b is an easy calculation.

For c, we consider Case IV (the other three being similar but easier).  We
note that we are in the situation where there are two peaks but the growth is
not maximal for the entire first half of the Hilbert function.

We know that
$h_{R/G'}(t) = h_{R/J}(t)$ for all $t \leq \ell$ (Lemma \ref{G hilb function
facts} (c)).  Hence  $G'$ and $J$ agree up degree
$\ell$.  Since $d_2 > \ell +1$, we get that the only generator is of degree
$d_1$ in this range.  (This holds for $G$, hence also for $G'$ by
semicontinuity.)  Hence we can compute the Hilbert function of
$R/G'$ as follows.  (We compute only the part that is relevant to our
subsequent calculation.)
\[
h_{R/G'}(t) = \left \{
\begin{array}{ll}
\binom{t+2}{2} & \hbox{if } 0 \leq t < d_1; \\ \\
\binom{t+2}{2} - \binom{t-d_1 +2}{2} & \hbox{if } d_1 \leq t \leq \ell; \\ \\
\hbox{(symmetric)} & \hbox{otherwise}
\end{array}
\right.
\]
It is known (cf.\ for instance \cite{diesel} Corollary 2.6) that  in each
degree, the third difference of the Hilbert function gives the smallest
possible number of minimal generators of $G'$ in that degree.  When the socle
degree is even, there is a $G'$ with precisely these minimal generators.  When
the socle degree $s$ is odd, one more generator may be needed, which occurs in
degree $\frac{s+3}{2}$, if the remaining number of generators is even (since
the total number of minimal generators must be odd -- cf.\ \cite{BE}).   In
Case IV, the socle degree is $d_1+d_2+d_3-d_4-3 = 2\ell +1$, which is odd. 
We compute
\[
\Delta^3 h_{R/G'}(t) = \left \{
\begin{array}{ll}
0 & \hbox{if } 0 < t < d_1 \\
-1 & \hbox{if } t = d_1 \\
-d_1 & \hbox{if } t = \ell + 1 \\
\geq 0 & \hbox{otherwise}
\end{array}
\right.
\]
This gives the result.
\end{proof}

Note that Case I is really a special case of Corollary \ref{one peak
resol}.  

\begin{theorem} \label{resol for codim 3}
Let $I = (G_1,G_2,G_3,G_4)$ be a general almost complete
intersection in $R$ of type $(d_1,d_2,d_3,d_4)$.  If $d_4 >
d_1+d_2+d_3-3$ then $G_4 \in (G_1,G_2,G_3)$ and $I$ is a complete
intersection, hence has a Koszul resolution.  So without loss of
generality assume that $d_1 \leq d_2 \leq d_3 \leq d_4 \leq
d_1+d_2+d_3-3$. Choose $F_1$ and $F_2$ to be the free modules
appearing in Proposition \ref{resol of general G} (depending on
the values of the $d_i$, i.e.\ depending on which of Cases I-IV
applies).  Then using these $F_i$, (\ref{codim 3 resol}) is a free
$R$-resolution for $I$.

Furthermore, the following represent all the splitting that occurs.
\begin{itemize}
\item[I.]  If $d_2+d_3 < d_1 +d_4 +3$, $d_1+d_2+d_3+d_4$ is odd {\em
and} $d_1+d_4 = d_2+d_3 -1$ then one summand
$R(d_1-d)$ splits at the end of (\ref{codim 3 resol}).

\item[II.] Assume that $d_2+d_3 < d_1+d_4 + 3$ and $d_1+d_2+d_3+d_4$ is even.
\begin{itemize}
\item[(i)] If $d_1+d_4 = d_2+d_3-2$ then one summand $R(d_1-d)$ splits at the
end of the resolution.

\item[(ii)] If $d_1 = d_2$, $d_3 = d_4$ and $\ell$ is even (where $\ell$ is
defined in Proposition \ref{resol of general G}) then one summand
$R(d_2-d)$ splits at the end of the resolution.

\item[(iii)] If $d_4-d_3 = d_2-d_1$ and $\ell$ is even then one summand
$R(d_1-d)$ splits at the end of the resolution.

\item[(iv)] If the hypothesis of (ii) holds then clearly the hypothesis of
(iii) also holds.  In this case there is only one splitting, by applying
either (ii) or (iii).  Other than this, no two of (i), (ii) and (iii) can
happen simultaneously.  
\end{itemize}

\item[III.] If $d_2+d_3 \geq d_1 +d_4 +3$ and $d_1+d_2+d_3+d_4$ is odd then one
summand $R(d_1-d)$ splits at the end of the resolution.

\item[IV.] If $d_2+d_3 \geq d_1 +d_4 +3$ and $d_1+d_2+d_3+d_4$ is even then one
summand $R(d_1-d)$ splits at the end of the resolution.

\end{itemize}
\end{theorem}

\begin{proof}
The approach is similar to that of Proposition \ref{aci two
peaks}.  In each of the four cases we start with the general
Gorenstein resolution described in Proposition \ref{resol of
general G} and link with general forms in $G$ of degrees
$d_1,d_2,d_3$ to produce an almost complete intersection $I$ of
type $(d_1,d_2,d_3,d_4)$.  Thanks to (\ref{codim 3 resol}) we have
a free resolution for $I$.  We want to produce from it the minimal
free resolution for a general almost complete intersection of the
same type.

Because the forms are general, if one of the $d_i$ is the degree of a minimal
generator of $G$, there is a corresponding splitting of the resolution.  (We
check below that the numerical conditions given in the statement of the
theorem give precisely this situation.)  This produces the claimed resolution
in the statement of the theorem.  It is possible that there are still summands
 that numerically are candidates for splitting, in this free resolution.

Suppose that the general almost complete intersection of this type has further
splitting.  Then linking back gives a Gorenstein ideal whose minimal free
resolution is smaller than the one described in Proposition \ref{resol of
general G}, which is a contradiction.

So we have only to check that the numerical conditions in the statement of the
theorem are exactly what is required to have a generator of degree $d_1$ or
$d_2$ (and that both cannot happen simultaneously).

We saw in Proposition \ref{resol of general G} (b) that $d_2 \geq \ell +2$, so
in cases I.\ and III.\ the only possibility is a splitting of a summand
$R(d_1-d)$.  This is automatic in III.\ and it is easy to check that in case
I.\ we have $d_1 = \ell +1$ if and only if $d_1 + d_4 = d_2 + d_3 -1$.

In case II.\ we can check that $d_1 = \ell +1$ if and only if $d_1 + d_4 = d_2
+ d_3 - 2$.  We can also check that $d_2 = \ell +2$ if and only if $d_2 + d_4 =
d_1 + d_3$, which in turn is equivalent to $d_1 = d_2$ and $d_3 = d_4$.
Clearly, then, we cannot have both $d_1 = \ell +1$ and $d_2 = \ell +2$.
Finally, we can also have $d_1 = \ell +2$ if $\ell$ is even, and this happens
if and only if $d_4 -d_3 = d_2 -d_1$.

In case IV.\ we automatically get a splitting of a summand $R(d_1-d)$.  We can
again check that $d_2 = \ell +2$ if and only if $d_1 = d_2$ and $d_3 = d_4$.
However, this time it is incompatible with the hypothesis $d_2+d_3 \geq d_1
+d_4 +3$.
\end{proof}

\begin{example} \label{codim 3 ghosts}
Combining Proposition \ref{resol of general G} and Theorem \ref{resol for codim
3} it is easy to give examples of general almost complete intersections which
have ``ghost'' terms in the resolution, i.e.\ summands which occur in
consecutive modules in the minimal free resolution but cannot be split off.
Some examples for
$(d_1,d_2,d_3,d_4)$ are $(4,4,4,8)$, $(5,5,6,8)$ and $(3,6,6,7)$.

When there is such a ghost term, it is a summand $R(\ell +2-d)$, which one
checks is $R(\frac{-d_1-d_2-d_3-d_4)}{2})$.

For instance, the minimal free $R$-resolution for $(4,4,4,8)$ is
\[
0 \rightarrow
\left (
\begin{array}{c}
R(-10) \\
\oplus \\
R(-11)^2
\end{array}
\right )
\rightarrow
\left (
\begin{array}{c}
R(-8)^3 \\
\oplus \\
R(-9)^2 \\
\oplus \\
R(-10)
\end{array}
\right )
\rightarrow
\left (
\begin{array}{c}
R(-4)^3 \\
\oplus \\
R(-8)
\end{array}
\right )
\rightarrow R \rightarrow R/I \rightarrow 0
\]
(here $\ell = 0$).  The summand $R(-8)$ in the first and second free modules is
clearly not going to split since one represents a minimal generator and the
other represents a Koszul relation between two other generators.  But the
summand
$R(-10)$ in the second and third free modules is an illustration of the
phenomenon that we are describing.  To our knowledge this is the first
counterexample to Iarrbino's Thin Resolution Conjecture \cite{iarrobino}.
\end{example}

A special case of interest is when the generators of our general almost
complete intersection all have the same degree, and here we record the minimal
free resolution in this case.

\begin{corollary} \label{aci same deg in 3 vars}
Let $I$ be a general almost complete intersection of type
$(a,a,a,a)$.  Then $R/I$ has a minimal free $R$-resolution
\[
0 \rightarrow R(-2a-1)^a \rightarrow
\left (
\begin{array}{c}
R(-2a)^3 \\
\oplus \\
R(-2a+1)^a
\end{array}
\right )
\rightarrow R(-a)^4 \rightarrow R \rightarrow R/I \rightarrow 0.
\]
\end{corollary}

\begin{proof}
The argument is a simple application of the preceding work.  The only thing
to note is that if $a$ is even then $\delta = 1$, and a summand $R(-2a)$ is
split off, but using a summand from $K_1^\vee (-d)$ rather than from $F_2^\vee
(-d)$.
\end{proof}

%%%%%%%%%%%%%%%%%%%%%%%%%%%%%%%%%%%%%%%%%%%%%%%%%%%%%%%%%%%%%%%%%%%%%%%%%%%

\section{Consequences of a more careful analysis of $R/G$}

In this section we will first extend Corollary \ref{aci same deg
in 3 vars} to rings of higher dimension.  Our strongest result comes when
this dimension $n$ is even and $a$ is also even, but we also have results for
the more general situation.  These are contained in Theorem \ref{main result
of section 5}.  These resolutions are very explicit.  Our methods apply more
generally, but the computations and notation become very cumbersome.  As a
middle road, we conclude the section with a result for $n=4$ which carefully
shows how to apply the method of this section and an inductive approach to
obtain the free resolution for $R/I$, which is minimal ``half'' of the time.

The following result from \cite{MN3} is crucial
to our work in this section:

\begin{proposition}[\cite{MN3} Proposition 8.7] \label{MN3 result}
Let $A = R/G$ be a graded Artinian Gorenstein $k$-algebra.  Let $L \in R$ be
a general linear form and let $\bar R = R/(L)$.  Let $s$ be the socle degree
of $A$, and let $\alpha = \init [0:_A L]$.  Then we have for all $i \in
{\mathbb Z}$
\[
\begin{array}{l}
\left [ \tor_i^R (A,k) \right ]_j = \\ \\

\hskip 1cm
\left \{
\begin{array}{ll}
\left [ \tor_i^{\bar R} (A/L A,k) \right ]_j
& \hbox{if } j \leq \alpha +i-2 \\
\leq \left [ \tor_i^{\bar R} (A/L A,k) \right ]_j +
\left [ \tor_{n-i}^{\bar R} (A/L A,k) \right ]_{s+n-j}
& \hbox{if } \alpha +i-1 \leq j \leq s-\alpha +i+1 \\
\left [ \tor_{n-i}^{\bar R} (A/L A,k) \right ]_{s+n-j}
& \hbox{if } j \geq s-\alpha +i+2
\end{array}
\right.
\end{array}
\]
\end{proposition}

This can be applied to our situation.  Since the calculations are somewhat
complicated, we illustrate it with two examples before we proceed to a more
general statement.

\begin{example} \label{case a=4}
Let $R = k[x_1,x_2,x_3,x_4]$ and let $I$ be a general almost complete
intersection of type $(4,4,4,4,4)$.  Let $J$ be the ideal given by the first
four generators and let $G = [J:I]$ be the linked Gorenstein ideal.  The
Hilbert function of $A := R/G$ is
\[
1 \ \ 4 \ \ 10 \ \ 20 \ \ 31 \ \ 20 \ \ 10 \ \ 4 \ \ 1 \ \ 0.
\]
 Let $L$ be a general linear form.  We have seen (Corollary \ref{wlp}) that $A
= R/G$ has the Strong Lefschetz property, hence the Weak Lefschetz property. 
So the Hilbert function of $A/LA$ is
\[
1 \ \ 3 \ \ 6 \ \ 10 \ \ 11 \ \ 0.
\]
We want to find the minimal free resolution of $A/LA$ over $\bar R =
R/(L) \cong k[x,y,z]$.

Note that the regularity of $A/LA$ is 5.  Hence it has a minimal
free $\bar R$-resolution of the form
\[
0 \rightarrow
\left (
\begin{array}{c}
\bar R(-6)^{c_1} \\
\oplus \\
\bar R(-7)^{c_2}
\end{array}
\right )
\rightarrow
\left (
\begin{array}{c}
\bar R(-5)^{b_1} \\
\oplus \\
\bar R(-6)^{b_2}
\end{array}
\right )
\rightarrow
\left (
\begin{array}{c}
\bar R(-4)^{a_1} \\
\oplus \\
\bar R(-5)^{a_2}
\end{array}
\right )
\rightarrow
\bar R \rightarrow A/LA \rightarrow 0.
\]
Clearly $a_1 = 4$, from the Hilbert function.  The four generators
in degree 4 are the restriction, $\bar J$, of $J$, hence can be
viewed as giving a general almost complete intersection in $\bar
R$.  $b_1$ represents all linear syzygies of these four
polynomials.   From Corollary \ref{aci same deg in 3 vars} we see
that $\bar J$ has no linear syzygies,  so $b_1 = 0$.  This  means
that the four generators in degree 4 span a subspace of dimension
$3 \cdot 4 = 12$ in $\bar R_5$, so we need $21 - 12 = 9$ more
generators in degree 5 and $a_2 = 9$.  Once $b_1 = 0$, this forces
$c_1 = 0$ since the smallest number is strictly increasing from
one free module to the next.  Then by exactness one computes $b_2
= 23$ and $c_2 = 11$, so the minimal free $\bar R$-resolution has
the form
\[
0 \rightarrow
\bar R(-7)^{11}
\rightarrow
\bar R(-6)^{23}
\rightarrow
\left (
\begin{array}{c}
\bar R(-4)^{4} \\
\oplus \\
\bar R(-5)^{9}
\end{array}
\right )
\rightarrow
\bar R \rightarrow A/LA \rightarrow 0.
\]

Now we can use Proposition \ref{MN3 result}.  In this case $s = 8$ and thanks
to the Weak Lefschetz property, $\alpha = 4$.  We have discussed $A/LA$ above
and have computed the minimal free $\bar R$-resolution of $A/LA$.  We conclude
(after a calculation) that $R/G$ has a minimal free resolution of the form
\[
0 \rightarrow R(-12) \rightarrow
\left (
\begin{array}{c}
R(-7)^{\leq 20} \\
\oplus \\
R(-8)^{\leq 4}
\end{array}
\right )
\rightarrow R(-6)^{\leq 46} \rightarrow
\left (
\begin{array}{c}
R(-4)^{\leq 4} \\
\oplus \\
R(-5)^{\leq 20}
\end{array}
\right )
\rightarrow R \rightarrow R/G \rightarrow 0.
\]
A computation from the Hilbert function gives that the inequalities are in
fact equalities.  Therefore we have the precise minimal free resolution of
$R/G$.

To compute the minimal free resolution of $R/I$, note that we have linked
using four quartics, which are minimal generators of $G$.  Therefore, using
(\ref{resol of R/I}) we see that we can split off four summands $R(-12)$ from
the end of the resolution of $R/I$.  The remaining resolution is
\[
0 \rightarrow R(-11)^{20} \rightarrow R(-10)^{46} \rightarrow
\left (
\begin{array}{c}
R(-8)^{10} \\
\oplus \\
R(-9)^{20}
\end{array}
\right )
\rightarrow R(-4)^5 \rightarrow R \rightarrow R/I \rightarrow 0
\]
which is minimal.
\end{example}

\begin{example} \label{case a=5}
Let $R = k[x_1,x_2,x_3,x_4]$ and let $I$ be a general almost complete
intersection of type $(5,5,5,5,5)$.  Let $J$ be the ideal given by the first
four generators and let $G = [J:I]$ be the linked Gorenstein ideal.  The
Hilbert function of $A := R/G$ is
\[
1 \ \ 4 \ \ 10 \ \ 20 \ \ 35 \ \ 52 \ \ 52 \ \ 35 \ \ 20 \ \ 10 \ \ 4 \ \ 1 \
\ 0.
\]
 Let $L$ be a general linear form.  We have seen (Corollary \ref{wlp}) that $A
= R/G$ has the Strong Lefschetz property, hence the Weak Lefschetz property. 
So the Hilbert function of $A/LA$ is
\[
1 \ \ 3 \ \ 6 \ \ 10 \ \ 15 \ \ 17 \ \ 0.
\]

We want to find the minimal free resolution of $A/LA$ over $\bar R =
R/(L) \cong k[x,y,z]$. Note that the regularity of $A/LA$ is 6.  Then precisely
the same reasoning as in Example \ref{case a=4} now gives that $A/LA$ has a
minimal $\bar R$-resolution of the form
\[
0 \rightarrow
\bar R(-8)^{17}
\rightarrow
\bar R(-7)^{36}
\rightarrow
\left (
\begin{array}{c}
\bar R(-5)^{4} \\
\oplus \\
\bar R(-6)^{16}
\end{array}
\right )
\rightarrow
\bar R \rightarrow A/LA \rightarrow 0.
\]

Now we can use Proposition \ref{MN3 result}.  In this case $s =
11$ and thanks to the Weak Lefschetz property, $\alpha  = 6$.  We
have discussed $A/LA$ above and have computed the minimal free
$\bar R$-resolution of $A/LA$.  We conclude (after a calculation)
that $R/G$ has a minimal free $R$-resolution of the form
\[
0 \rightarrow R(-15) \rightarrow
\left (
\begin{array}{c}
R(-8)^{y} \\
\oplus \\
R(-9)^{\leq 16} \\
\oplus \\
R(-10)^4
\end{array}
\right )
\rightarrow
\left (
\begin{array}{c}
R(-7)^{19+y} \\
\oplus \\
R(-8)^{19+y}
\end{array}
\right )
 \rightarrow
\left (
\begin{array}{c}
R(-5)^{4} \\
\oplus \\
R(-6)^{\leq 16} \\
\oplus \\
R(-7)^y
\end{array}
\right )
\rightarrow R \rightarrow R/G \rightarrow 0
\]
where $y \leq 17$.  A computation from the Hilbert function gives that the
inequality ``$\leq 16$'' is in fact an equality, but it does not determine
the value of $y$.  Experiments with Macaulay \cite{macaulay} indicate that
in fact $y=0$, but we have not been able to prove this.

To compute the minimal free resolution of $R/I$, note that we have linked
using four quartics, which are minimal generators of $G$.  Therefore, using
(\ref{resol of R/I}) we see that we can split off four summands $R(-12)$ from
the end of the resolution of $R/I$.  The remaining resolution is
\[
0 \rightarrow
\left (
\begin{array}{c}
R(-14)^{16} \\
\oplus \\
R(-13)^y
\end{array}
\right )
 \rightarrow
\left (
\begin{array}{c}
R(-13)^{19+y} \\
\oplus \\
R(-12)^{19+y}
\end{array}
\right )
 \rightarrow
\left (
\begin{array}{c}
R(-10)^{10} \\
\oplus \\
R(-11)^{16} \\
\oplus \\
R(-12)^y
\end{array}
\right )
\rightarrow R(-5)^5 \rightarrow R \rightarrow R/I \rightarrow 0.
\]
As mentioned above, computations with Macaulay \cite{macaulay} indicate that
$y=0$.
\end{example}

This approach leads us to the following more general results, for the minimal
free resolution of a general almost complete intersection.  We first consider
the case of general forms of the same degree, $a$.
 Our  goal is to extend Corollary \ref{aci same deg in 3 vars} to rings of
higher dimension. To this end we need to introduce some extra notation.
For any integers $0<a , n \in \ZZ$, we set
\[
\begin{array}{rcl}
s(n,a) & = & (n-1)a-n \\ \\ \ell (n,a) & = & \displaystyle \left
\lfloor \frac{s(n,a)}{2} \right \rfloor \\ \\ t(n,a) & = & \max
\{t \mid \ell (n,a)+t-1\ge ta \}
\\
\\ \alpha _j(n,a) & = & \displaystyle \sum _{i=0}^{t(n,a)} (-1)^{i+j-1}
\binom{n}{i}\binom{\ell (n,a)+n-2+j-ia}{n-2}-\sum _{r=1}^{j-1}
 (-1)^{r+j}
\binom{n-2+j-r }{j-r} \alpha _r(n,a) \\ \\ &  &  \mbox{ defined
inductively for } 1\le j \le n-2
\end{array}
\]
\[
\begin{array}{rcl}
\\ \alpha_{n-1} & = & \displaystyle  \sum _{i=2}^{n-t(n,a)-1}(-1)^{i}
\alpha _{n-i}(n,a)+\sum _{i=n-t(n,a)}^{n-1} (-1)^{i}(\alpha
_{n-i}(n,a) +\beta _{n-i}(n))+(-1)^{n}
\\
\\ \beta_{i}(n) & = &\displaystyle \binom{n}{i} \mbox{ for } i\ge
1
\\ \\ \displaystyle \binom{a}{b} & = & 0 \mbox{ if } a<b.
\end{array}
\]

\begin{theorem} \label{main result of section 5}
Let $I$ be a general almost complete intersection in $R =
k[x_1, \dots ,x_n]$ of type $(a, \dots ,a)$, $a >1$.

If $n$ is odd  then $R/I$ has a minimal free $R$-resolution
\[
\begin{array}{l}
0 \rightarrow \left (
\begin{array}{c}
R(\ell (n,a)+1-na) ^{\alpha _1(n,a)} \\ \oplus \\ R(\ell
(n,a)+2-na) ^{\le \alpha _{n-1}(n,a)}
\end{array}
\right )
 \rightarrow  \dots 
 \rightarrow
\left (
\begin{array}{c}
  R(\ell (n,a)+(\frac{n-1}{2})-na) ^{\le \alpha _{\frac{n-1}{2}}(n,a)}
   \\ \oplus \\
R(\ell (n,a)+(\frac{n+1}{2})-na) ^{\le \alpha
_{\frac{n+1}{2}}(n,a) }
\end{array}
\right ) \rightarrow
\end{array}
\]
\[
\begin{array}{l}
 \left (
\begin{array}{c} R(-\frac{n+1}{2} a)^{\beta_{\frac{n+1}{2}}(n)} \\
\oplus \\
  R(\ell (n,a)+(\frac{n+1}{2})-na) ^{\le \alpha
  _{\frac{n+1}{2}}(n,a)}
   \\ \oplus \\
R(\ell (n,a)+(\frac{n+3}{2})-na) ^{\le \alpha
_{\frac{n-1}{2}}(n,a)}
\end{array}
\right ) \rightarrow \left (
\begin{array}{c} R(-\frac{n-1}{2} a)^{\beta_{\frac{n-1}{2}}(n+1)} \\
\oplus \\
  R(\ell (n,a)+(\frac{n+3}{2})-na) ^{\le \alpha
  _{\frac{n+3}{2}}(n,a)}
   \\ \oplus \\
R(\ell (n,a)+(\frac{n+5}{2})-na) ^{\le \alpha
_{\frac{n-3}{2}}(n,a)}
\end{array}
\right ) \rightarrow  \dots  \rightarrow
\end{array}
\]
\[
\begin{array}{l}
 \left (
\begin{array}{c} R(-2a)^{\beta_{2}(n+1)} \\
\oplus \\
  R(\ell (n,a)+n-1-na) ^{\le \alpha _{n-1}(n,a)}
   \\ \oplus \\
R(\ell (n,a)+n-na) ^{ \alpha _{1}(n,a)}
\end{array}
\right ) \rightarrow
 R(-a)^{\beta_{1}(n+1)} \rightarrow R \rightarrow R/I \rightarrow 0
\end{array}
\]

\vskip 4mm  If $n$ is even and $s(n,a)$ is odd  then $R/I$ has a
minimal free $R$-resolution

\[
\begin{array}{ll}
\begin{array}{ll}
 0 & \rightarrow \left (
\begin{array}{c}
R(\ell (n,a)+1-na) ^{\alpha _1(n,a)} \\ \oplus \\ R(\ell
(n,a)+2-na) ^{\le \alpha _{n-1}(n,a)}
\end{array}
\right )
 \rightarrow \dots  \\ \\
& \rightarrow
\left (
\begin{array}{c}
  R(\ell (n,a)+(\frac{n}{2})-na) ^{\le \alpha _{\frac{n}{2}}(n,a)}
   \\ \oplus \\
R(\ell (n,a)+(\frac{n}{2})+1-na) ^{\le \alpha _{\frac{n}{2}}(n,a)}
\end{array}
\right ) \rightarrow
\end{array}
\\ \\
\begin{array}{ll}
& \rightarrow \left (
\begin{array}{c} R(-\frac{n}{2} a)^{\beta_{\frac{n}{2}}(n+1)} \\
\oplus \\
  R(\ell (n,a)+(\frac{n}{2})+1-na) ^{\le \alpha
  _{\frac{n}{2}+1}(n,a)}
   \\ \oplus \\
R(\ell (n,a)+(\frac{n}{2})+2-na) ^{\le \alpha
_{\frac{n}{2}-1}(n,a)}
\end{array}
\right ) \rightarrow \dots
\end{array}
\\ \\
\begin{array}{ll}
& \rightarrow
 \left (
\begin{array}{c} R(-2a)^{\beta_{2}(n+1)} \\
\oplus \\
  R(\ell (n,a)+n-1-na) ^{\le \alpha _{n-1}(n,a)}
   \\ \oplus \\
R(\ell (n,a)+n-na) ^{ \alpha _{1}(n,a)}
\end{array}
\right )
 \rightarrow
 R(-a)^{\beta_{1}(n+1)} \rightarrow R \rightarrow R/I \rightarrow 0
\end{array}
\end{array}
\]

\vskip 2mm If $n$ is even and $s(n,a)$ is even   then $R/I$ has a
minimal free $R$-resolution
\[
\begin{array}{ll}
0 & \rightarrow R(\ell (n,a)+1-na) ^{\alpha
_1(n,a)+\alpha_{n-1}(n,a)}
 \rightarrow \dots
 \rightarrow
  R(\ell (n,a)+(\frac{n}{2})-na) ^{2 \alpha _{\frac{n}{2}}(n,a)}
\\ \\
& \rightarrow
\left (
\begin{array}{c} R(-\frac{n}{2} a)^{\beta_{\frac{n}{2}}(n+1)} \\
\oplus \\
  R(\ell (n,a)+(\frac{n}{2})+1-na) ^{ \alpha
  _{\frac{n}{2}+1}(n,a)
   + \alpha _{\frac{n}{2}-1}(n,a)}
\end{array}
\right ) \rightarrow \dots \\ \\
& \rightarrow
 \left (
\begin{array}{c} R(-2a)^{\beta_{2}(n+1)} \\
\oplus \\
  R(\ell (n,a)+n-1-na) ^{\alpha _{n-1}(n,a)+
 \alpha _{1}(n,a)}
\end{array}
\right )
 \rightarrow
 R(-a)^{\beta_{1}(n+1)} \rightarrow R \rightarrow R/I \rightarrow 0
\end{array}
\]
\end{theorem}

\begin{proof} We proceed by induction on $n$. The case  $n=3$ is
covered by Corollary \ref{aci same deg in 3 vars} and, even more, all possible
splitting does occur.

For arbitrary $n$, one first checks that the complete intersection
$J$ formed by taking $n$ of the generators of $I$ links $I$ to a
Gorenstein ideal $G$ with socle degree $s(n,a)$. According to
Lemma \ref{G hilb function facts}, if $s(n,a)$ is even then there
is only one peak, occurring in degree $\ell (n,a) $.  If $s(n,a)$
is odd then there are two peaks, the first of which is in degree
$\ell (n,a)$. Let $A = R/G$ and let $L \in R_1$ be a general
linear form.

\vskip 2mm The Hilbert function of $R/G$ is
\[
h_{R/G}(t) = \left \{
\begin{array}{ll}
\binom{t+n-1}{n-1} & \hbox{if } t \leq a-1 \\ \\
\binom{t+n-1}{n-1} - n \cdot \binom{t-a+n-1}{n-1} & \hbox{if } a
\leq t \leq \ell (n,a)\\ \\ (symmetric) & \hbox{otherwise }.
\end{array}
\right.
\]
By the Weak Lefschetz property, the Hilbert function of $A/LA$ is
\[
h_{A/LA} (t) = \left \{
\begin{array}{ll}
\binom{t+n-2}{n-2} & \hbox{if } t \leq a-1 \\ \\
\binom{t+n-2}{n-2} - n \cdot \binom{t-a+n-2}{n-2} & \hbox{if } a
\leq t \leq \ell (n,a)\\ \\ 0 & \hbox{if } t > \ell (n,a).
\end{array}
\right.
\]
Note that if $\bar R = R/(L)$ and $\bar G = \frac{G+(L)}{(L)}$
then $A/LA \cong \bar R / \bar G$.  We will compute the minimal
free $\bar R$-resolution of $A/LA$.

We first observe that the regularity of $A/LA$ is $\ell (n,a)+ 1 =
\lfloor \frac{(n-1)a-n}{2} \rfloor +1$. Note also that $(A/LA)_i
\cong (\bar R / \bar G)_i = (\bar R / \bar J)_i$ for all $i \leq
\ell (n,a)$ by Lemma \ref{G hilb function facts}.  Since $\bar
R/\bar J$ is an Artinian almost complete intersection in $\bar R$,
by hypothesis of induction we have very good bounds on the graded
Betti numbers of $\bar R/\bar J$.

We have seen that the minimal generators of $G$ in degree $\leq
\ell (n,a)$ agree with those in $J$ of degree $\leq \ell (n,a)$,
hence all appear in degree $a$. Furthermore, since $\reg A/LA =
\ell (n,a)+1$, the minimal free $\bar R$-resolution of $A/LA$
begins
\[
\dots \rightarrow
\begin{array}{c}
\bar R (-a)^n \\ \oplus\\ \bar R(-\ell (n,a)-1)^{\alpha_1(n,a)}
\end{array}
\rightarrow \bar R \rightarrow A/LA \rightarrow 0.
\]

\vskip 2mm {\em Claim:} (a) If $1\le i\le t(n,a)$ then
\[
[\tor_i^{{\bar R}}(A/LA,k)]_{j} = \left \{
\begin{array}{ll}
\beta_{i}(n) & \hbox{if } j = ia \\  \alpha _{i}(n,a)  & \hbox{if
} j = \ell (n,a) +i  \\ 0 & \hbox{otherwise.}\\
\end{array}
\right.
\]
(b) If $n-1\ge i>t(n,a)$ then
\[
[\tor_i^{{\bar R}}(A/LA,k)]_{j} = \left \{
\begin{array}{ll}
  \alpha _{i}(n,a)  & \hbox{if } j = \ell (n,a) +i
\\ 0 & \hbox{otherwise.}\\
\end{array}
\right.
\]

\vskip 2mm {\em Proof of the Claim:}

The generators of $A/LA$ occur in degrees $a$ and $\ell (n,a)+1$
and  the generators of degree $a$ of $A/LA$ are precisely the
generators of $\bar J$. So
 $(A/LA)_t \cong   (\bar R /
\bar J)_t$ for all $t \leq \ell (n,a)$ and we deduce that
$[\tor_i^{\bar R} (A/LA,k)]_j=[\tor_i^{\bar R} (\bar R/\bar
J,k)]_j$ for all $j\le \ell (n,a)+i-1$. Furthermore,
 we have very good bounds on the graded
Betti numbers of $\bar R/\bar J$ by the inductive hypothesis.  In
particular,
\[
[\tor_i^{\bar R} (\bar R/\bar J,k)]_j \neq 0 \Leftrightarrow \left
\{
\begin{array}{l}
j = ia \hbox{ and } i \leq \lfloor \frac{n+1}{2} \rfloor \\ 
j = - \ell(n-1,a) - (n-1) + (n-1)a + i - 1 \\ 
j = - \ell(n-1,a) - (n-1) + (n-1)a + i - 2 \hbox{ (possibly)}
\end{array}
\right.
\]
and,
\[
[\tor_i^{\bar R} (\bar R/\bar J,k)]_{ia}=\beta_{ia}(n) \mbox{ if }
i \leq \lfloor \frac{n-1}{2} \rfloor .
\] 
An easy calculation shows
that
\[
 \ell(n,a) +i-1< -\ell(n-1,a) - (n-1) + (n-1)a + i-2
\]
and \[  \ell(n,a) +i-1<ia \Leftrightarrow i>t(n,a)
\]
which together with the fact that
 $reg (A/LA) = \ell(n,a)+1$ give us that the values
of $j$ given in part (a)  and (b) of the claim are the only ones
where $[\tor_i^{\bar R} (A/LA,k)]_j \neq 0$.  The precise values
of $\beta_i(n)$ are given by the free ${\bar R}$-resolution of
$\bar R/\bar J$ and the precise values of  $\alpha _{i}(n,a)$ then
follow from a calculation involving the Hilbert function of $A/AL$
and the exactness of the $\bar R$-resolution.

\vskip 2mm It follows from the above claim that $A/LA$ has a
minimal free ${\bar R}$-resolution of the form
\[
0 \rightarrow \bar R(-\ell (n,a) -n+1)^{\alpha_{n-1}(n,a)}
\rightarrow  \dots \rightarrow \bar R(-\ell (n,a)
-t(n,a)-1)^{\alpha_{t(n,a)+1}(n,a)} \rightarrow
\]
\[
\begin{array}{c}
\bar R (-t(n,a)a)^{\beta_{t(n,a)}(n)} \\ \oplus \\ \bar R(-\ell
(n,a)-t(n,a))^{\alpha _{t(n,a)}(n,a)}
\end{array}
\rightarrow  \dots  \rightarrow
\begin{array}{c}
\bar R (-a)^{\beta_1(n)} \\ \oplus \\ \bar R(-\ell
(n,a)-1)^{\alpha _1(n,a)}
\end{array}
\rightarrow \bar R \rightarrow A/LA \rightarrow 0.
\]

Now we take the now-known resolution for $A/LA$ and plug the
values into Proposition \ref{MN3 result}. When $n$ is odd we obtain (the other
two cases are similar), after a calculation, that

If $1\le i\le t(n,a)$ then
\[
[\tor_i^{ R}(A,k)]_{j} = \left \{
\begin{array}{ll}
\beta_{i}(n) & \hbox{if } j = ia \\ \le  \alpha _{i}(n,a)  &
\hbox{if } j = \ell (n,a) +i \\ \le  \alpha _{n-i}(n,a)  &
\hbox{if } j = \ell (n,a) +i+1
 \\ 0 & \hbox{otherwise,}\\
\end{array}
\right.
\]

\noindent  if $\frac{n-1}{2}\ge i>t(n,a)$ then
\[
[\tor_i^{R}(A,k)]_{j} = \left \{
\begin{array}{ll}
 \le \alpha _{i}(n,a)  & \hbox{if } j = \ell (n,a) +i \\
  \le \alpha _{n-i}(n,a)  & \hbox{if } j = \ell (n,a) +i+1
\\ 0 & \hbox{otherwise}\\
\end{array}
\right.
\]

\noindent and, by symmetry (since $G$ is Gorenstein) we compute
$[\tor_i^{R}(A,k)]_{j} $ for $\frac{n+1}{2}\le i \le n-1$.

From the Hilbert function and the above considerations we can
compute that the number of minimal generators of $G$ in degree
$\ell (n,a) +1$ is precisely $\alpha_1$. So we have equality for
$[\tor_1^{R}(A,k)]_{\ell (n,a) +1}$. But then by symmetry (since
$G$ is Gorenstein) we get equality for $[\tor_{n-1}^{
R}(A,k)]_{\ell (n,a) +n}$.

Finally, (\ref{resol of R/I}) gives the claimed free
$R$-resolution of $R/I$ after splitting off the summands of the
left half part of the resolution since they correspond to minimal
generators of $J$ and syzygies involving these minimal generators.
When $n$ and $s(n,a)$ are both even this resolution is clearly
minimal.
\end{proof}

\begin{example}
Suppose $n=4$.  Notice that the case $a=2$ is covered by Corollary
\ref{one peak resol}, the case $a=3$ by Corollary \ref{codim 4 result}, the
case $a=4$ by Example \ref{case a=4} and the case $a=5$ by Example \ref{case
a=5}.  So suppose that $a \geq 6$ and assume that $a$ is even.  Then
\[
\begin{array}{rcl}
s=s(4,a) & = & 3a-4 \\ \\ 
\ell=\ell(4,a) & = & \displaystyle
\frac{3a-4}{2}  \\ \\ 
t=t(n,a) & = & 1 \\ \\ 
\alpha_1=\alpha_1(4,a) & = & \displaystyle 
\binom{\ell+3}{2} - 4 \binom{\ell+3-a}{2} \\ \\ 
\alpha_2=\alpha_2(4,a) & = &
\displaystyle - \binom{\ell +4}{2} + 4  \binom{\ell +4-a}{2} + 3
\alpha_1 \\ \\ 
\alpha_3=\alpha_3(n,a) & = &\displaystyle \alpha _2 - \alpha_1 -
3
\end{array}
\]
and $R/I$ has a minimal free $R$-resolution
\[
\begin{array}{l}
\displaystyle 0 \rightarrow R \left (\ell+1-4a \right )^{\alpha_
1+\alpha_ 3} \rightarrow R \left (\ell+2-4a \right )^{2\alpha _2}
\rightarrow \left (
\begin{array}{c}
R(-2a)^{10} \\ \oplus \\ \displaystyle R \left ( \ell+3-4a \right
)^{\alpha _1+\alpha _3}
\end{array}
\right ) \hbox{\hskip 2cm}
\\ \\
\hfill \rightarrow R(-a)^5 \rightarrow R \rightarrow R/I
\rightarrow 0.
\end{array}
\]
\end{example}

It should be clear from the above considerations that for any values of $n$
and of $d_1,\dots,d_{n+1}$ we can say quite a bit about the minimal free
resolution using our methods, and that in some cases we can give the precise
resolution.  This is of course the most desirable, and interesting,
situation.  It is also clear that the notation gets progressively more
cumbersome as $n$ grows.  As a middle road, we show that the precise minimal
free resolution can be obtained when $n=4$ and when we have only one peak,
otherwise allowing the values of the $d_i$ to be arbitrary.  

\begin{theorem} \label{general ans when n=4}
Assume that $n=4$.  Let $I = (G_1,\dots, G_5)$ be a general almost
complete intersection of type $(d_1,\dots,d_5)$ and assume that
$\sum_{i=1}^5 d_i$ is even, with $2 \leq d_1 \leq \dots \leq d_5
\leq \sum_{i=1}^4 d_i - 4$. (The latter condition only assures
that $I$ is not a complete intersection.)  Then the minimal free
resolution of $R/I$ can be obtained by the above procedure, and is
explicitly given in the proof below.
\end{theorem}

\begin{proof}
Let $J = (G_1,\dots,G_4)$ and let $G = [J:I]$ as usual.  We have
$d = \sum_{i=1}^4 d_i$, $s = d-d_5-4$ and $\ell = \frac{s}{2}$.
Note that we are in the case of one peak, which occurs in degree
$\ell$.  We have seen in Lemma 2.6 that $h_{R/G}(t) = h_{R/J}(t)$
for $t \leq \ell$.

Let $L$ be a general linear form and let $A = R/G$.  For a
numerical function $f$, we again denote by $\Delta f$ the first
difference function $ \Delta f = f(t) - f(t-1)$ for $t \in
{\mathbb Z}$. By the Weak Lefschetz property for $A$ (Corollary
2.7), we get
\[
h_{A/LA}(t) = \left \{
\begin{array}{ll}
\Delta h_{R/J}(t) & \hbox{for } t \leq \ell; \\ 0 & \hbox{for } t
> \ell.
\end{array}
\right.
\]
In particular the regularity of $A/LA$ is $\ell +1$.

Let $\bar R = R/(L)$ and $\bar J = \frac{J+(L)}{(L)}$.   Let $f =
d_1+d_2+d_3$.  If $d_4 \leq d_1 + d_2 + d_3 -3$ then a free $\bar
R$-resolution for $\bar R / \bar J$ is given by 
\begin{equation} \label{resol of bar J}
0 \rightarrow F_1^\vee (-f) \rightarrow \left (
\begin{array}{c}
\bigoplus_{i=1}^3 \bar R(d_i -f) \\ \oplus \\ F_2^\vee (-f)
\end{array}
\right ) \rightarrow \bigoplus_{i=1}^4 \bar R(-d_i) \rightarrow
\bar R \rightarrow \bar R/\bar J \rightarrow 0
\end{equation}
(with all splitting as described in Theorem 4.2).  If $d_4 > d_1 + d_2 +
d_3 -3$ then $\bar J$ is a complete intersection; this case is easier and we
leave it to the reader.

Our first task is to describe the minimal free ${\bar R}$-resolution of
$A/LA$.  Since $G_t = J_t$
for $t \leq \ell$, we have that $\bar R /\bar J \cong A/LA$ in
degrees $\leq \ell$, and $h_{A/LA}(t) = h_{\bar R/\bar J}(t)$ for
$t \leq \ell$.  Since the regularity of $A/LA$ is $\ell +1$, we
have
\begin{equation} \label{tor of A/LA}
[\tor_i^{\bar R} (A/LA,k)]_j = \left \{
\begin{array}{ll}
[\tor_i^{\bar R} (\bar R/\bar J,k)]_j , &  j \leq \ell +i-1 \\ ? &
j = \ell +i \\ 0 & j > \ell +i
\end{array}
\right.
\end{equation}
We have seen that there may be overlap (``ghost terms'') in the
minimal free resolution of $\bar R / \bar J$, but these have been
completely described in terms of the $d_i$.  They are all
contained in the first line of (\ref{tor of A/LA}).  No further
overlaps can arise from the second line.

Since $A/LA$ ends in degree $\ell$, we have in particular that 
\[
\begin{array}{rcl}
[\tor_3^{\bar R}(A/LA,k)]_{\ell +3} & = & h_{A/LA}(\ell) \\ & = &
\Delta h_{R/J}(\ell):=b_3
\end{array}
\]
Using Proposition \ref{resol of general G}, Theorem \ref{resol for codim 3} 
and the resolution (\ref{resol of bar J}), it is tedious but possible to
check that
\[
[\tor_3^{\bar R}(\bar R/\bar J,k)]_j = 0 \ \hbox{ for } j \leq
\ell +2.
\]
Combining the above, we obtain a precise description of
$[\tor_3^{\bar R}(A/LA,k)]_j$ for all $j$.

Again using Proposition \ref{resol of general G}, Theorem \ref{resol for
codim 3}  and the resolution (\ref{resol of bar J}), we get that if
$j\le \ell +1$, then there is only one possibility for $[\tor_2^{\bar R}
(\bar R/\bar J,k)]_j \neq 0$:
\[
\begin{array}{rcl}
[\tor_2^{\bar R} (\bar R/\bar J,k)]_j \neq 0 & \Leftrightarrow 
&  j = f-d_3 \leq \ell +1 \\
& \Leftrightarrow & d_1 + d_2 + d_5 + 2 \leq d_3 + d_4.
\end{array}
\]
Notice that even if $f-d_3 = \ell +1$, there is an overlap in the resolution
which does not split.  This comes from Theorem \ref{resol for codim 3}. 
Such a summand corresponds to a non-trivial syzygy of generators of lower
degree.

Therefore $A/LA$ has the following minimal free ${\bar R}$-resolution
\begin{equation}\label{resol of A/LA}
0 \rightarrow {\bar R}(-\ell -3)^{b_3} \rightarrow \left (
\begin{array}{c}
 \bar R(d_3 -f)^{a_{2}} \\ \oplus \\ {\bar
R}(-\ell -2)^{b_2}
\end{array}
\right ) 
\rightarrow \left (
\begin{array}{c}
\bigoplus_{i=1}^4 \bar R(-d_i)^{a^i_{1}} \\ \oplus \\ {\bar
R}(-\ell -1)^{b_1} 
\end{array} 
\right ) 
\rightarrow \bar R \rightarrow A/LA \rightarrow 0 
\end{equation} 
where 
\[
a_{2}=\left
\{
\begin{array}{l}
1 \hbox{ if } f-d_3 \leq \ell +1 \\ 
0 \hbox { otherwise }
\end{array}
\right.
, \hskip 1cm a^i_{1}=\left \{
\begin{array}{l}
1 \hbox{ if } d_i \leq \ell
\\ 0 \hbox { otherwise, } 
\end{array}
\right. 
\]
and $b_2$,  $b_1$ are determined by the exactness of the
above exact sequence (use the equation involving only ranks and
the equation involving the first Chern classes).  Note that $b_3$ was
defined above.

 Now we take the now-known ${\bar R}$-resolution of $A/LA$ and
plug the values into   Proposition 5.1. Note that $\alpha$ is what
we now are calling $\ell$ (because there is only one peak), and in
fact $\alpha = \ell = s-\alpha$.  Then, we have
\[
\begin{array}{l}
\left [ \tor_i^R (A,k) \right ]_j = \\ \\

\hskip 1cm \left \{
\begin{array}{ll}
\left [ \tor_i^{\bar R} (A/L A,k) \right ]_j & \hbox{if } j \leq
\ell +i-2 \\ \leq \left [ \tor_i^{\bar R} (A/L A,k) \right ]_j +
\left [ \tor_{4-i}^{\bar R} (A/L A,k) \right ]_{2\ell+4-j} &
\hbox{if } \ell +i-1 \leq j \leq \ell +i+1 \\ \left [
\tor_{4-i}^{\bar R} (A/L A,k) \right ]_{2\ell+4-j} & \hbox{if } j
\geq \ell +i+2.
\end{array}
\right.
\end{array}
\]
Now, using this fact together with the fact that $A$ and $R/J$ agree in
degree $\leq \ell$, plus the symmetry of the resolution, we get that the
minimal free $R$-resolution of $R/G$ has the form
\[
0 \rightarrow R(-2\ell -4) \rightarrow  \left (
\begin{array}{c}
\bigoplus_{i=1}^4  R(d_i-2\ell -4)^{a_1^{i}} \\ \oplus \\
 R(-\ell -3)^{\le b_1+b_3}
\end{array} 
\right )
\rightarrow 
\left (
\begin{array}{c}
R(d_3 -f)^{a_2} \\ \oplus \\  R(-\ell -2)^{\le 2b_2} \\ \oplus \\
 R(f-d_3 -2 \ell-4)^{a_2}
\end{array}
\right )
\]

\[\rightarrow \left (
\begin{array}{c}
\bigoplus_{i=1}^4  R(-d_i)^{ a_1^{i}} \\ \oplus \\  R(-\ell
-1)^{\le b_1+b_3} \end{array} \right ) \rightarrow  R \rightarrow
R/G \rightarrow 0 \]
(We are also again using the fact that if $d_3 -f = \ell +1$ then this
syzygy represents a minimal syzygy of forms of lower degree, and cannot be
split off.)

We now claim that the above inequalities are in fact equalities.  Indeed, we
know that $A$ and $R/J$ agree up to degree $\ell$.  Let us denote by $J_{\leq
\ell}$ the ideal generated by the components of $J$ in degree $\leq \ell$. 
Note that $J_{\leq \ell}$ is a complete intersection of height $<4$, hence
$\hbox{depth } R/J_{\leq \ell} \geq 1$.  It follows that we need $b_1$
minimal generators in degree $\ell +1$ in order to bring $\dim (R/J_{\leq
\ell})_{\ell +1}$ down to the level of $\dim A_\ell$, and then another $b_3$
generators to bring it down to the level of $\dim A_{\ell +1}$ ($= \dim
A_{\ell -1}$).  Therefore we have the first and third free modules.  The
equality for the second free module is exactly obtained from the ranks of
the minimal free resolution (\ref{resol of A/LA}).  Therefore we have the
precise minimal free resolution of $R/G$.

To compute the minimal free $R$-resolution of $R/I$, note that we
have linked using $G_1$, $G_2$, $G_3$ and $G_4$ and that $G_i$ is
a minimal generator of $G$ if  $d_i\le \ell $ or if $d_i = \ell +1$.
Therefore, using (\ref{resol of R/I}) we see that for any $i$ such that
$d_i\le \ell +1 $ we can split off a summand $R(d_i-d)$ from the end of the
resolution of $R/I$ and the remaining resolution is minimal.  The fact that
no further splitting occurs follows as in the first part of the proof of
Theorem \ref{resol for codim 3}.
\end{proof}

\begin{example} \label{3 3 4 6 6}
Let $(d_1,\dots,d_5) = (3,3,4,6,6)$.  Then the above proof gives a minimal
free resolution for $R/G$ to be
\[
0 \rightarrow R(-10) \rightarrow 
\left (
\begin{array}{c}
R(-6)^{17} \\
\oplus \\
R(-7)^2
\end{array}
\right )
\rightarrow R(-5)^{36} \rightarrow 
\left (
\begin{array}{c}
R(-3)^2 \\
\oplus \\
R(-4)^{17} 
\end{array}
\right )
\rightarrow R \rightarrow R/G \rightarrow 0
\]
and a minimal free resolution for $R/I$ to be
\[
0 \rightarrow R(-12)^{16} \rightarrow 
\left (
\begin{array}{c} 
R(-10) \\
\oplus \\
R(-11)^{36}
\end{array}
\right )
\rightarrow 
\left (
\begin{array}{c}
R(-6) \\
\oplus \\
R(-7)^2 \\
\oplus \\
R(-9)^4 \\
\oplus \\
R(-10)^{18}
\end{array}
\right )
\rightarrow
\left (
\begin{array}{c}
R(-3)^2 \\
\oplus \\
R(-4) \\
\oplus \\
R(-6)^2
\end{array}
\right )
\rightarrow R \rightarrow R/I \rightarrow 0.
\]
Notice that all the summands $R(-6), R(-7)$ and $R(-9)$ are Koszul, but that
there is also an $R(-10)$ that does not split.  
\end{example}

We observe that all of the examples where redundant (``ghost'') terms occurred
were cases where the degrees were not all the same.  In contrast, Corollary
\ref{aci same deg in 3 vars} and Theorem \ref{main result of section 5} at least
suggest the following conjecture (and prove some cases):

\begin{conjecture}
Let $I \subset R = k[x_1,\dots,x_n]$ be the ideal of $n+1$ generically chosen
forms of the same degree.  Then there is no redundant term in the minimal free
resolution of $R/I$.  Consequently, the minimal free resolution can be computed
from the Hilbert function, whose conjectured value is given by Fr\"oberg.
\end{conjecture}

\noindent Note that the first syzygy case of this conjecture  is proved by
Hochster and Laksov
\cite{hochster-laksov}.

%\end{enumerate}


\begin{thebibliography}{999}

\bibitem{anick} D. Anick, {\em Thin Algebras of embedding dimension three},
J.\ Algebra {\bf 100} (1986), 235--259.

\bibitem{aubry} M. Aubry, {\em S\'erie de Hilbert d'une alg\`ebre
de polynomes quotient}, J.\ Algebra {\bf 176} (1995), 392--416.

\bibitem{macaulay} D.\ Bayer and M.\ Stillman, Macaulay: A system for
computation in algebraic geometry and commutative algebra. Source and object
code available for Unix and Macintosh computers.  Contact the authors, or
download from ftp://math.harvard.edu via anonymous ftp.

\bibitem{bigatti} A.\ Bigatti, {\em Upper bounds for the Betti numbers of a
given Hilbert function}, Comm.\ Algebra {\bf 21} (1993), no. 7, 2317--2334.

\bibitem{boij} M.\ Boij, {\em Betti numbers of compressed level algebras}, J.\
Pure and Applied Algebra {\bf134} (1999), 111--131.

\bibitem{boij2} M.\ Boij, {\em Gorenstein Artin Algebras and Points in
Projective Space}, Bull.\ London Math.\ Soc. {\bf 31} (1999), 11--16.

\bibitem{C} K.\ Chandler, {\em The Geometric Fr\"oberg-Iarrobino
Conjecture}, in preparation.

\bibitem{BE} D.\ Buchsbaum and D.\ Eisenbud, {\em Algebra Structures for Finite
Free Resolutions, and some Structure Theorems for Ideals of Codimension 3},
Amer.\ J.\ of Math.\ {\bf 99} (1977), 447--485.

\bibitem{DGO} E.\ Davis, A.V.\ Geramita, F.\ Orecchia, {\em Gorenstein Algebras
and the Cayley-Bacharach Theorem}, Proc.\ Amer.\ Math.\ Soc.\ {\bf 93} (1985),
593--597.

\bibitem{diesel} S.~Diesel, {\em Irreducibility and Dimension Theorems for
Families of Height 3 Gorenstein Algebras}, Pacific J.\ Math. {\bf 172}
(1996), no. 2, 365--397.

\bibitem{EP} D.\ Eisenbud and S.\ Popescu, {\em Gale Duality and Free Resolutions of
Ideals of Points}, Invent.\ math.\ {\bf 136} (1999), 419--449.

\bibitem{froberg} R. Fr\"oberg, {\em An inequality for Hilbert series of
graded algebras}, Math.\ Scand.\ {\bf 56} (1985), 117--144.

\bibitem{froberg-hollman} R. Fr\"oberg, J. Hollman, {\em Hilbert series for
ideals generated by generic forms}, J.\ Symbolic Comput.\ {\bf 17}
(1994), 149--157.

\bibitem{harima} T.\ Harima, {\em Characterization of Hilbert functions of
Gorenstein Artin algebras with the weak Stanley property}, Proc.\ Amer.\
Math.\ Soc.\ {\bf 123} (1995), 3631--3638.

\bibitem{HMNW} T.\ Harima, J.\ Migliore, U.\ Nagel and J.\ Watanabe, {\em The
Weak and Strong Lefschetz properties for Artinian $K$-algebras}, preprint.

\bibitem{HS} A.\ Hirschowitz and C.\ Simpson, {\em La r\`esolution minimale
de l'id\`eal d'un arrangement g\`en\`eral d'un grand nombre de points dans
$\proj{n}$},Invent. \ Math. {\bf 126} (1996), 467--503.

\bibitem{hochster-laksov} M.\ Hochster and D.\ Laksov, {\em The Linear
Syzygies of Generic  Forms}, Comm.\ Algebra {\bf 15} (1987), 227--239.

\bibitem{hulett} H.\ Hulett, {\em Maximum Betti numbers of homogeneous ideals
with a given Hilbert function}, Comm.\ Algebra {\bf 21} (1993), no. 7,
2335--2350.

\bibitem{iarrobino} A. Iarrobino, {\em Inverse system of a symbolic power III.
Thin algebras and fat points}, Compos. Math. {\bf 108} (1997) 319-356.

\bibitem{IK} A.\ Iarrobino and V.\ Kanev, ``Power Sums, Gorenstein
Algebras, and Determinantal Loci,'' Springer LNM 1721 (1999).

\bibitem{lorenzini} A.\ Lorenzini, {\em The Minimal Resolution Conjecture},
J.\ Alg.\ {\bf 156} (1993), 5--35.

\bibitem{migliore} J.\ Migliore, ``Introduction to Liaison Theory and
Deficiency Modules,''  Birkh\"auser, Progress in Mathematics 165, 1998.

\bibitem{MN3} J.\ Migliore and U.\ Nagel, {\em Reduced arithmetically
Gorenstein schemes and simplicial polytopes with maximal Betti numbers},
preprint.

\bibitem{PR} Keith Pardue and Ben Richert, {\em Resolutions of Generic
Ideals}, in preparation.

\bibitem{PS} C.\ Peskine and L.\ Szpiro, {\em Liaison des vari\'et\'es
alg\'ebriques.\ I}, Inv.\ Math.\ {\bf 26} (1974), 271--302.

\bibitem{stanley} R.\ Stanley, {\em Weyl groups, the hard Lefschetz
theorem, and the Sperner property}, SIAM J.\ Algebraic Discrete Methods
{\bf 1} (1980), 168--184.

\bibitem{watanabe} J.\ Watanabe, {\em The Dilworth number of Artinian
rings and finite posets with rank function}, Commutative Algebra and
Combinatorics, Advanced Studies in Pure Math.\ Vol.\ 11, Kinokuniya Co.\
North Holland, Amsterdam (1987), 303--312.

\end{thebibliography}
\end{document}